\documentclass[11pt,leqno]{amsart}
\usepackage[dvips]{pstcol}
\usepackage{enumerate,amsthm,amsmath,amssymb,comment,graphicx,psfrag, float}
\usepackage{xspace}
\usepackage{mathrsfs}
\usepackage{xcolor}
\usepackage{fancybox}
\usepackage{subfigure}
\usepackage{bm}
\usepackage{enumerate, verbatim}
\input amssym.def
\numberwithin{equation}{section}

\begin{document}
\newcommand{\Ind}[1]{\mbox{ \bf 1}_{#1}} 
\newcommand{\Dem}{\noindent {\it {Proof.} \ }}
\newcommand{\BE}{\begin{equation}}
\newcommand{\dds}[1]{\frac{d {#1}}{ds}}
    \renewcommand{\topfraction}{0.9}	
    \renewcommand{\bottomfraction}{0.8}	
    \setcounter{topnumber}{2}
    \setcounter{bottomnumber}{2}
    \setcounter{totalnumber}{4}     
    \setcounter{dbltopnumber}{2}    
    \renewcommand{\dbltopfraction}{0.9}	
    \renewcommand{\textfraction}{0.07}	
    \renewcommand{\floatpagefraction}{0.7}	
    \renewcommand{\dblfloatpagefraction}{0.7}	


\def\F{{\cal{F}}}
\newcommand{\E}{{\mathcal{E}}}
\def\<{{\langle }}
\def\>{{\rangle }}
\def\V{{\mathcal{V}}}
\def\sg{{\varSigma_\gamma}}
\def\Th{{\mathcal{T} }_h }
\def\Sh{{\mathcal{S} }_h }
\def\Nh{{\mathcal{N} }_h }

\def\H{ { L^2 (\varOmega) } }
\def\V{ { H^1_0 (\varOmega) } }

\def\Po{{P^1_h}}
\def\Pz{{P^0_h}}

\def\al{\alpha}
\def\wU{\widehat U}
\def\wA{\widehat A}
\def\wa{\widehat a}
\def\gam{\gamma}
\def\be{\beta}
\def\la{\lambda}
\def\th{\vartheta}
\def\tth{\tilde{\vartheta}}
\def\eps{\epsilon}
\def\ds{\displaystyle}
\def\Co{{\Bbb C}}
\def\Re {{\Bbb R}}
\def\P{{\Bbb P}}
\def\N{{\Bbb N}}
\def\R+{{\Bbb R}_*^+}
\def\Ze{{\Bbb Z}}
\def\Na{{\Bbb N}}
\def\bV{{\Bbb V}}

\newcommand{\rec}{{\mathcal{R}}}
\def\tdelta{{\tilde\varDelta}}

\def\tm{\tilde{m}}
\def\tp{\tilde{p}}

\def\proof{\smallskip
\noindent
{\em Proof.\ } }

\def\Lemma{\smallskip
\noindent
{\bf Proof.\ } }

\addtolength{\subfigcapskip}{-0.1cm}

\newtheorem{lemma}{Lemma}[section]
\newtheorem{theorem}{Theorem}[section]
\newtheorem{remark}[lemma]{Remark}
\newtheorem{corollary}{Corollary}[section]
\newtheorem{notation}[lemma]{Notation}
\newtheorem{definition}{Definition}[section]
\newtheorem{subsec} {} [section]
\vspace*{2.0cm}      \renewcommand{\aa}{{\mbox{\boldmath$a$}}}
\newtheorem{assumption}{Assumption}[section]

\newcommand{\bb}{{\mbox{\boldmath$b$}}}

\newcommand{\nn}{{\mbox{\boldmath$n$}}}
\newcommand{\pp}{{\mbox{\boldmath$p$}}}
\newcommand{\snn}{{\mbox{\boldmath$\scriptstyle n$}}}

\newcommand{\xxi}{{\mbox{\boldmath$ \xi$}}}

\newcommand{\sxxi}{{\mbox{\boldmath$\scriptstyle \xi$}}}

\newcommand{\zzeta}{{\mbox{\boldmath$ \zeta$}}}

\newcommand{\eeta}{{\mbox{\boldmath$\eta$}}}

\newcommand{\da}{\downarrow}

\newcommand{\ua}{\uparrow}

\newcommand{\Proof}{{\sc Proof.}\ }

\def\Box{\vbox{\offinterlineskip\hrule

        \hbox{\vrule\phantom{o}\vrule}\hrule}}

\newcommand{\bull}{\rule{2mm}{2mm}}

\newcommand{\var}{\mathop{\rm Var}\nolimits}

\renewcommand{\div}{\mathop{\rm div}\nolimits}
\newcommand{\sign}{\mathop{\rm sign}\nolimits}
\newcommand{\argmin}{\mathop{\rm arg\,min}\nolimits}

\newcommand{\jacfield}{{\mathfrak{A}}}

\renewcommand{\a}{\mathfrak{a}}

\renewcommand{\eps}{\varepsilon}

\newcommand{\interfunct}{{\Phi}}

\newcommand{\der}{\delta}

\newcommand{\gield}{{\mathfrak{G}}}

\newcommand{\dder}{{\der U}}

\newcommand{\fdata}{F}                  

\newcommand{\coerc}{\varSigma}             

\newcommand{\Coerc}{\varrho}               

\newcommand{\coercgamma}{\coerc_\gamma}  

\newcommand{\ttau}{{\pc\tau}}

\newcommand{\field}{{\mathfrak{F}}}

\newcommand{\vfield}{{\field_{\V}}}

\newcommand{\aield}{{\mathfrak{A}}}

\newcommand{\bield}{\mathfrak{B}}

\newcommand{\iield}{\mathfrak{I}}

\newcommand{\resolvent}{{\mathfrak{J}}}

\newcommand{\residual}{{\cal R}}

\newcommand{\sig}{{{\boldsymbol{\varSigma}}}}
\renewcommand{\H}{{{\bf H}}}
\newcommand{\A}{{{\bf A}}}
\renewcommand{\V}{{{\bf V}}}
\renewcommand{\u}{{{\bf u}}}
\renewcommand{\v}{{{\bf v}}}
\newcommand{\n}{{{\bf n}}}
\newcommand{\w}{{{\bf w}}}
\newcommand{\Z}{{{\bf Z}}}
\newcommand{\bP}{{{\bf P}}}
\newcommand{\bR}{{{\bf R}}}
\newcommand{\J}{{{\bf J}}}
\newcommand{\f}{{{\bf f}}}
\newcommand{\U}{{{\bf U}}}
\newcommand{\W}{{{\bf W}}}
\newcommand{\e}{{{\bf e}}}
\newcommand{\h}{{{\bf h}}}
\newcommand{\z}{{{\bf z}}}
\newcommand{\g}{{{\bf g}}}
\newcommand{\bepsilon}{{\boldsymbol {\epsilon}}}
\newcommand{\brho}{{\boldsymbol{\rho}}}

\newcommand{\bphi}{{\boldsymbol{\varphi}}}
\newcommand{\bPhi}{{\boldsymbol{\Phi}}}
\newcommand{\bchi}{{\boldsymbol{\chi}}}

\newcommand{\Vs}{{{\cal V^*}}}

\newcommand{\Vsz}{{{\cal V}^*_0}}

\newcommand{\partition}{{\cal P}}       

\newcommand{\pc}[1]{{\bar{#1}}}         

\newcommand{\ppc}[1]{\mbox{\kern-3pt\b{\kern3pt$#1$}}}


\newcommand{\pl}[1]{{#1}}

\newcommand{\estuno}{{\cal E}}

\newcommand{\estdue}{{\cal D}}

\newcommand{\estduet}{{\tilde\estdue}}

\newcommand{\Err}{E}

\newcommand{\Errga}{\Err_\gamma}

\newcommand{\err}{e}

\newcommand{\TAU}{\tau}

\newcommand{\lip}{L}

\newcommand{\dQ}{{\partial\varOmega\times(0,T)}}

\newcommand{\Q}{{\bf Q}}

\renewenvironment{description}[1][0pt]
  {\list{}{\labelwidth=0pt \leftmargin=#1
   \let\makelabel\descriptionlabel}}
  {\endlist}
\title[A posteriori estimates for fully discrete fractional-step $\vartheta$-approximations]
{A posteriori error estimates for fully discrete fractional-step $\vartheta$-approximations
for the heat equation}

\author[\sf F. Karakatsani]{Fotini \ Karakatsani}
\address{Department of  Mathematics and Statistics,
26 Richmond street, G1 1XH,
Glasgow, UK
} 
\email {\sf{foteini.karakatsani{\it @\,}strath.ac.uk} }

\begin{abstract}
We derive optimal order a posteriori error estimates for 
fully discrete approximations of  the initial-boundary value problem for the heat equation. For the discretization in  time 
we apply the fractional-step $\th$-scheme and for the discretization in space 
the finite element method with finite element spaces that are allowed to change 
with time.  The first optimal order a posteriori error estimates  in $L^\infty(0,T;L^2(\varOmega))$ are derived by
applying the reconstruction technique.

\end{abstract}



\maketitle

\section{Introduction}
Let $\varOmega$ be a  bounded polyedral domain in $\Re^d, d = 2, 3,$ with boundary $\partial\varOmega$ and $T > 0.$ In the present paper we derive optimal order a posteriori error estimates in $L^\infty(0,T;L^2(\varOmega))$  for fully discrete fractional-step $\th$-scheme approximations  for the heat equation:
\begin{equation}
\label{heat_eq}
\left \{
\begin{alignedat}{2}
&u_t -\varDelta u = f  \quad && \text{ \rm{in} } \varOmega \times (0,T], \\
&u=0  \quad &&\text{ \rm{on} } \partial\varOmega \times (0,T], \\
&u(x,0) = u^0(x) \quad &&\text{ \rm{in} } \varOmega.\\
\end{alignedat}
\right .
\end{equation}
We assume throughout that $f\in L^2(0,T; L^2(\varOmega))$ and $u^0 \in H_0^1(\varOmega);$   then the weak solution  $u$ of \eqref{heat_eq} belongs to $L^\infty(0,T; H_0^1(\varOmega))$  with 
$\partial_t u\in L^2(0,T; L^2(\varOmega)).$

Adaptive finite element methods  are a fundamental numerical tool  in computational science and engineering for approximating partial differential equations  with solutions that exhibit nontrivial characteristics.  They aim to automatically adjust the mesh to fit the numerical solution, that means  fine meshes in the regions where the solution changes fast and coarse in the regions where the solution changes slowly and, consequently, to keep the computational cost as low as possible.  The general structure of a loop of an adaptive algorithm for evolution equations is: 
Given the approximation $U^{n}\in V_{n}$  (reflecting the space discretization method)
\begin{description}\label{basic_step_AFEM}
\item[\rm{1\ }] choose the next time node $t^{n+1}$ and the next space $V_{n+1},$ 
  \item[\rm{2\ }] project $U^{n}$  to  $V_{n+1}$ to get $\tilde U^{n},$ 
 \item[\rm{3\ }] use  $\tilde U^{n}$ as starting value to perform the evolution step in  $V_{n+1}$ to obtain the
new approximation $U^{n+1} \in V_{n+1}\, .$
\end{description}
The design of such algorithms, particularly the  decision made at the first step of the loop,  
is usually based on suitable a posteriori error estimates which can measure the quality of 
the approximate solution and provide information of the error distribution.  

Although the fractional-step $\th$-scheme was first proposed  as an  operator splitting method 
in the context of time-dependent Navier--Stokes equations, cf.\ \cite{Glowinski85}, 
\cite{Glowinski86}, and \cite{Bri-Glo_Per87}, it is an attractive alternative to popular 
time-stepping schemes, cf.\ \cite{Rannacher98}, 
 \cite{Glowinski:2003}.  Indeed, its parameters can be  chosen such that  to produce a 
strongly A-stable and  second order accurate method.  Thus, the scheme can combine the second-order  accuracy
of the Crank--Nicolson method with the full smoothing property of the backward Euler method
in the case of non-smooth initial data.
Moreover, in contrast to the backward Euler, it is very little
numerically dissipative and, compared to the Runge--Kutta methods of higher
order, of lower complexity and storage requirements. 
 For more details  we refer to \cite{Urbaniak93}, \cite{Mueller94}, \cite{Rannacher98}, 
 \cite{Glowinski:2003} and the references therein. 

Despite the great deal of effort that has been devoted to the a posteriori error analysis of linear or nonlinear parabolic equations,  cf., for example, \cite{JNT}, \cite{EJ_l1}, 
\cite{EJ_nonl}, \cite{NSaV}, \cite{MN2003}, \cite{Verfuerth2003}, \cite{AMN2006}, \cite{LM2006}, 
\cite{BKM2012}, the results in case of  the  fractional-step $\th$-scheme are limited, to our 
knowledge at least, cf.  \cite{K2012}, \cite{MR2014}. 
Particularly,   a posteriori error estimates of optimal order in $L^\infty(0,T; L^2(\varOmega))$ 
were derived for time-discrete approximations of linear parabolic equation in \cite{K2012}.  
The key for the a posteriori error analysis was the use of 
a continuous piecewise quadratic in time approximation of $u,$  the so-called \emph{fractional-step $\th$-reconstruction}, whose residual was 
second order accurate. The definition of the fore-mentioned reconstruction followed 
the idea of the  \emph{two-point Crank-Nicolson reconstruction}, cf., \cite{AMN2006}.

Here, following the ideas developed  in \cite{BKM2012, BKM2013}, we combine  the fractional-step $\vartheta$-scheme and the Galerkin finite
element method to get a fully discrete 
scheme consistent with the mesh modification. The first optimal order a posteriori error estimate in $L^\infty(0,T; L^2(\varOmega))$ for fully fractional-step $\vartheta$-approximations are derived by exploiting both the elliptic reconstruction, cf.  \cite{MN2003},  and time-reconstruction techniques, cf.  \cite{AMN2006}, \cite{K2012} and \cite{LPP2009}. In particular,  we define a continuous representation  $\hat\omega,$ $\hat\omega:[0,T]\rightarrow H_0^1(\varOmega),$  of the approximate solution $U$  which will be referred  to as \emph{space-time reconstruction} of $U.$ The space-time reconstruction $\hat\omega$ is a piecewise quadratic polynomial in time which is based on approximations on either one time subinterval 
or two adjacent time subintervals.    
Then, the \emph{total error} $e:=u-U$ may be  split as 
\begin{equation*}
e=u-U = (u- \hat \omega) +(\hat \omega - U) =:
\hat\rho + \varepsilon,
\end{equation*}
where 
\begin{enumerate}[\small $\bullet$] \itemsep = 0pt
\item the \emph{space-time reconstruction error} $\varepsilon$ may be split as the sum of the  \emph{elliptic reconstruction error} and the \emph{time reconstruction error}. The elliptic reconstruction error can be bounded by using any elliptic estimator at our disposal and the time reconstruction error can be controlled by a posteriori quantities of optimal order.
\item the \emph{parabolic error} $\hat\rho$ satisfies an appropriate heat equation whose right hand-side can be bounded by computable quantities of optimal order.
\end{enumerate}

The rest of the paper is organized as follows. In 
Section 2 we introduce notation and the fully discrete scheme allowing mesh change. In
Section 3 we first discuss the space- and time-discretization and the corresponding reconstructions and then we present the space-time reconstruction. Specific choices of the reconstructions 
leading to estimators based on approximations on one time subinterval and on two adjacent time subintervals are given. Section 4 is devoted to the error analysis of the parabolic error $\hat\rho$ and  we state the final estimates in both aforementioned cases  of time-space reconstructions.  
The  asymptotic behavior of the derived estimators is presented in  Section 5.

\section{Preliminaries}

In this section we introduce the necessary notation for our analysis and the fully discrete scheme.

\subsection{Notation}
Let $0=t^0< t^1< \cdots <t^N=T,$ $I_n:=(t^{n-1},t^n]$ and $k_n:=t^n-t^{n-1}.$
For $\th\in(0, \frac 13),$ $\tth=1-2\th,$ we 
introduce the intermediate time levels $t^{n-1+\th} = t^{n-1} +\th k_n$ and 
$t^{n-\th} = t^{n-1} + (\th+\tth) k_n.$ 

For each $0\leq n \leq N,$  let  $\mathcal{T}_n$ be a triangulation of $\varOmega$ into disjoint $d$-simplices $K$ and   $h_n$ its {\it local mesh-size function}   defined by
\begin{equation}
\label{def_hn}
h_n(x):=\text{diam}(K), \quad K \in \mathcal{T}_n \text{ and } x\in K.
\end{equation} 
We assume that  the aspect ratios of all the elements are uniformly bounded with respect to $n=0, \ldots, N,$ and  the intersection of two different elements is either empty, or consists of a common vertex,  a common
edge, or  a common face.

We associate with 
each $\mathcal{T}_n$ a finite element space $\bV_h^n$  
\begin{equation}
\widetilde{\Bbb{V}}_h^n:=\{\phi\in H^1(\varOmega): \forall K\in\mathcal{T}_n: \phi|_K\in \Bbb{P}^l\} \quad \text{and} \quad {\Bbb{V}}_h^n:=\widetilde{\Bbb{V}}_h^n\cap H_0^1(\varOmega),
\end{equation}
where $\Bbb{P}^l$ is the space of polynomials in $d$ variables of degree at 
most $l.$  

For each $n$ and for each $K\in \mathcal{T}_n,$ we denote by 
$\mathcal{E}_n(K)$ the set of the facets of $K$ 
and  by $\varSigma_n(K)\subset \mathcal{E}_n(K)$ 
the set of the interior facets of $K,$ that is the facets that do not belong to the boundary of $\varOmega.$  
In addition, we introduce the sets 
$\mathcal{E}_n:=\cup_{K\in\mathcal{T}_n}\mathcal{E}_n(K)$ and
 $\varSigma_n:=\cup_{K\in\mathcal{T}_n}\varSigma_n(K)$.  

We also make the assumption that all triangulations are derived from the same macro-triangulation by using an admissible refinement procedure, e.g., the bisection-based refinement procedure used in ALBERTA-FEM toolbox, cf.\ \cite{SS2005}.  
Given two successive triangulations  $\mathcal{T}_{n-1}$ and 
$\mathcal{T}_n,$ we define the {\it finest common coarsening} 
$\hat{\mathcal{T}}_n:= \mathcal{T}_{n-1}\wedge\mathcal{T}_n$  and  the {\it coarsest common refinement} $\check{\mathcal{T}}_n:= \mathcal{T}_{n-1}\vee \mathcal{T}_n,$ whose local mesh-sizes are 
given by $\hat{h}_n:=\max(h_n,h_{n-1})$ and $\check{h}_n:=\min(h_n,h_{n-1}),$ respectively.  Note that essentially $\hat{\mathcal{T}}_n$ is the triangulation of ${\Bbb{V}}_h^n\cap\Bbb{V}_h^{n-1} $ and $\check{\mathcal{T}}_n$ is the triangulation of ${\Bbb{V}}_h^n\cup{\Bbb{V}}_h^{n-1}.$  In addition, we shall denote by $\check\varSigma_n$ and $\hat{\varSigma}_n$ the sets of the interior facets which correspond to $\check{\mathcal{T}}_n$ and  $\hat{\mathcal{T}}_n,$ respectively, namely $\check{\varSigma}_n:=\cup_{K\in \check{\mathcal{T}}_n}\varSigma_n(K)$ and $\hat{\varSigma}_n:=\cup_{K\in \hat{\mathcal{T}}_n}\varSigma_n(K).$ We refer to  \cite{LM2006}
for precise definitions.

We shall use
the shorthand notation $u^m(\cdot) :=u(\cdot,t^m)$ and $f^m(\cdot):=f(\cdot,t^m)$ throughout.  The jump  $J[\bm v]_e$  of a  discontinuous vector valued function $\bm v$ across an interior facet $e\in\varSigma_n(K)$  is defined by 
\begin{equation}
J[\bm v]_e (x):= \lim_{\delta\rightarrow 0} [\bm v(x+\delta \bm n_e) - \bm v(x-\delta \bm n_e)]\cdot \bm n_e,\end{equation}
where $\bm n_e$ is a unit normal vector  on $e$ and $x\in e.$

We denote by $\<\cdot, \cdot\>$  either the inner product in  $L^2(\varOmega)$ or the 
duality pairing  between $H^{-1}(\varOmega)$ and  $H_0^1(\varOmega)$  and  we let $a(\cdot,\cdot)$ be defined as $a(u, v):=\<\nabla u, \nabla v\>.$ 
For $\mathcal{D}\subset \mathbb{R}^d$ we denote 
by 
$\|\cdot\|_\mathcal{D}$ the 
norm in $L^2(\mathcal{D}),$ by $\|\cdot\|_{r,\mathcal{D}}$ 
and by $|\cdot|_{r,\mathcal{D}} $ the norm and the semi-norm,
respectively, in the Sobolev space $H^r(\mathcal{D}), r\in\mathbb{N}.$ 
In view of the Poincar\'e inequality, we consider $|\cdot|_{1,\mathcal{D}}$ to be the norm
in $H_0^1(\mathcal{D})$ and denote by  $|\cdot|_{-1,\mathcal{D}}$ 
the norm in $H^{-1}(\mathcal{D});$ whenever $\mathcal{D}=\varOmega$ the subscript $\mathcal{D}$ will be omitted  in the notation of function spaces and  norms.  

In addition, we shall use  the following notation for  functions $v$ defined in a piecewise sense
\begin{equation}
\label{elementwise_norms}
\begin{aligned}
 \|h_n^i v\|_{\mathcal{T} _n}& = \left(\sum _ {K \in \mathcal{T} _n} \|h_K^i v\|_{K} ^2 \right)^{1/2} \\ \quad 
 \text{and} \quad
 \|h_n^{i+\frac 12}J[\nabla v]\|_{\varSigma_n} &= \left(\sum _ {e \in \varSigma_n} \|h_e^{i+\frac 12}J[\nabla v]_e\|^2_{e}\right)^{1/2}, \quad i=1,2.
\end{aligned}
\end{equation}
\subsection{Discrete operators and interpolants} For $0\leq n \leq N,$ let  $\varDelta_h^n: H_0^1 \rightarrow {\Bbb{V}}_h^n$ be the discrete Laplacian 
corresponding to the finite element space ${\Bbb{V}}_h^n$ defined by 
\begin{equation}
\label{dis_laplacian} 
\<(-\varDelta_h^n) v, \chi_n\> =\<\nabla v, \nabla \chi_n\> \quad \forall \chi_n \in {\Bbb{V}}_h^n.
\end{equation}
Moreover,  we denote by $P_0^n,$ $P_0^n : L^2 \rightarrow {\Bbb{V}}_h^n,$ the $L^2$-projection onto ${\Bbb{V}}_h^n$ and,  in order our analysis to include several possible choices for the projection step, $\varPi^n, \widetilde\varPi^n: \Bbb{V}_h^{n-1}  \rightarrow  \Bbb{V}_h^n $ will  denote appropriate projections or interpolants to be chosen.

We now recall  the stability property and the approximation properties of the Cl\'ement-type interpolant introduced in  \cite{scott-zhang}. 

\begin{lemma}\label{lem_clem_app1} Let $\mathcal{I}^n: H_0^1 \rightarrow \Bbb{V}_h^n$ be a Cl\'ement-type interpolant. Then, we have
\begin{equation}
\label{clement_stab}
|\mathcal{I}^n z|_1 \leq c_1 |z|_1.
\end{equation} 
Furthermore, for $j\leq l+1,$ the following approximation properties are satisfied 
\begin{equation}
\label{clement_app1}
\begin{aligned}
&\|h_n^{-j}(z-\mathcal{I}^nz)\|_{\mathcal{T}^n}\le c_{1,j}|z|_j, & \\
&\|h_n^{1/2-j}(z-\mathcal{I}^nz)\|_{\varSigma_n}\le c_{2,j}|z|_j ,
\end{aligned}
\end{equation}  
where $l$ is the finite element polynomial degree and the constants
$c_1,$ $c_{1,j}$ and $c_{2,j}$ depend only on the shape-regularity of the 
family of triangulations $\{\mathcal{T}_n\}_{n=0}^N.$ 
\qed
\end{lemma}

Let  $C_E$ denote the elliptic regularity constant, that is
\begin{equation}
\label{ell_reg}
|v| _2 \leq C_E \| \varDelta v\|, \quad v\in H^ 2 (\varOmega)\cap H^ 1_0 (\varOmega)\, ,
\end{equation}
and $c_1, $ $c_{i,j}, i=1,2,\; j\leq l+1,$ be the constants in Lemmas \ref{lem_clem_app1}.
We shall also
use the notation 
$$ C_{j, 2}: = C_E\,  c_{j, 2}\, $$ for the constants appeared in the definition of the a posteriori error estimators.

\subsection{The fully discrete scheme}

We discretize \eqref{heat_eq} by 
applying the following Galerkin fractional-step $\vartheta$-scheme (GFS-scheme):
for a given approximation $U^0$ of $u^0$ and for $1\leq n \leq N,$ find $U^n\in {\Bbb{V}}_h^n,$ such that
\begin{equation}
\label{theta_scheme}
\left \{
\begin{aligned}
&\frac{U^{n-1+\th}-\varPi^n U^{n-1}}{\th k_n} - \alpha_1\varDelta_h^n U^{n-1+\th}
-\beta_1\widetilde\varPi^n \varDelta_h^{n-1} U^{n-1}= 
P_0^n \{ \alpha_2  f^{n-1+\th} + \beta_2 f^{n-1}\},\\
&\frac{U^{n-\th}-U^{n-1+\th}}{\tth k_n} - \beta_1 \varDelta_h^n U^{n-\th}
-\alpha_1 \varDelta_h^n U^{n-1+\th}= P_0^n \{ \beta_2  f^{n-\th} 
+ \alpha_2 f^{n-1+\th}\}, \\
&\frac{U^n-U^{n-\th}}{\th k_n}  - \alpha_1\varDelta_h^n U^n -
\beta_1 \varDelta_h^n U^{n-\th}= P_0^n \{
\alpha_2  f^{n} + \beta_2   f^{n-\th} \},
\end{aligned}
\right .
\end{equation}
with $\alpha_1, \alpha_2 \in (0,1),$ and $\beta_1 = 1-\alpha_1, \beta_2 = 1-\alpha_2.$  We shall sometimes find it convenient to rewrite 
\eqref{theta_scheme} in the form
\begin{equation}\label{theta_scheme1}
\left \{
\begin{aligned}
&\frac{U^{n-1+\th}-\varPi^n U^{n-1}}{k_n}- \bigl\{\alpha_1 \th \varDelta_h^n U^{n-1+\th}
+\beta_1 \th\widetilde\varPi^n \varDelta_h^{n-1} U^{n-1}\bigr\} \\
 =& 
 P_0^n \bigl\{\alpha_2 \th  f^{n-1+\th} + \beta_2 \th  f^{n-1}\bigr\}, \\
&\frac{U^{n-\th}-\varPi^nU^{n-1}}{k_n} -\bigl\{\beta_1 \tth \varDelta_h^n U^{n-\th}
+\alpha_1 (\th+\tth) \varDelta_h^n U^{n-1+\th}+\beta_1 \th \widetilde\varPi^n \varDelta_h^{n-1} U^{n-1}\bigr\}
\\=  & P_0^n \bigl\{\beta_2 \tth  f^{n-\th}  +  \alpha_2 (\th+\tth) f^{n-1+\th} 
+ \beta_2 \th f^{n-1}\bigr\},\\
&\frac{U^n-\varPi^nU^{n-1}}{k_n}  - \bigl\{ \alpha_1 \th\varDelta_h^n U^n+ \beta_1 (\th+\tth) \varDelta_h^n U^{n-\th}
+\alpha_1 (\th+\tth) \varDelta_h^n U^{n-1+\th} \\ & +\beta_1 \th \widetilde\varPi^n \varDelta_h^{n-1} U^{n-1}\bigr\}
\\=&
  P_0^n \bigl\{\alpha_2 \th  f^{n}  +  \beta_2 (\th+\tth) f^{n-\th} + 
 \alpha_2 (\th+\tth) f^{n-1+\th} + \beta_2 \th  f^{n-1}\bigr\}.
\end{aligned}
\right .
\end{equation}

Throughout the rest of the paper we shall assume that $\th=1-\frac{\sqrt 2}{2}$ and $\alpha_1 \in (\frac 1 2, 1],$ which implies that the fractional-step $\th$-scheme  is second-order accurate and $A(0)$-stable.  Indeed, the assumption that 
$\alpha_1 \in (\frac{1}{2}, 1]$ implies the strong $A$-stability of our scheme.  Furthermore, we can easily seen that the quadrature rule 
\begin{equation}\label{quad1}
\mathcal{I}_{\alpha,\th}(\phi) :=\beta\th \phi(0) +
\alpha (\th+\tilde{\th})\phi(\th) + \beta(\th +\tilde{\th})\phi(1-\th) +
\alpha \th \phi(1)\approx \int_0^1 \phi (s)ds
\end{equation}
integrates first degree polynomials exactly if and only if
$\alpha=\beta=\frac 1 2$ or $\th=1-\frac{\sqrt{2}}{2}$.
Thus, the  assumption $\th=1-\frac{\sqrt 2}{2}$  ensures that 
the fractional-step $\th$-scheme 
is second-order accurate with respect to time.  We refer to \cite{Glowinski:2003} for more details.

\subsection{The fully discrete scheme in compact form} 
We introduce the following piecewise linear polynomials with respect to time 
\begin{equation}
\label{def_phi}
 \varphi(t):=\ell_0^n(t)f^{n-1} + \ell_1^n(t) f^{n}, \quad t\in I_n,
\end{equation}
and
\begin{equation}
\label{def_Theta}
\varTheta(t)=\ell_0^n(t)\widetilde\varPi^n(-\varDelta_h^{n-1})U^{n-1} + \ell_1^n(t) (-\varDelta_h^n)U^n, \quad t\in I_n,
\end{equation}
with 
\begin{equation}
\label{def_l0_l1}
\ell_0^n(t) := \frac{t^n-t}{k_n} \quad \text{and} \quad \ell_1^n(t) := \frac{t-t^{n-1}}{k_n},
\quad t\in I_n.
\end{equation}
Moreover, we let $\hat\varTheta$ be defined as 
\begin{equation}\label{def_hatTheta}
\begin{aligned}
\widehat{\varTheta}(t):= \varTheta(t) -  \xi_\varTheta^n, \quad t \in I_n,
\end{aligned}
\end{equation}
with 
\begin{equation}\label{corr_Theta}
\xi_\varTheta^n:=  (\th+\tth) \bigl\{\alpha
[\varTheta(t^{n-1+\th})+\varDelta_h^nU^{n-1+\th}]+\beta
[\varTheta(t^{n-\th})+\varDelta_h^nU^{n-\th})]\bigr\},\quad t \in I_n,
\end{equation}
and $\hat\varphi$
\begin{equation}\label{def_hatphi}
\hat{\varphi}(t):= \varphi(t) -  \xi_\varphi^n,\quad t \in I_n,
\end{equation}
with
\begin{equation}\label{corr_phi}
\xi_\varphi^n:=  (\th+\tth) \bigl\{\alpha
[\varphi(t^{n-1+\th})-f^{n-1+\th}] +\beta
[\varphi(t^{n-\th})-f^{n-\th})]\bigr\},\quad t \in I_n.
\end{equation}

Note that both $\xi_\varTheta^n$ and $\xi_\varphi^n$ are  a posteriori quantities of optimal order, cf. \cite{K2012} for details. 

According to definitions \eqref{def_hatTheta} and \eqref{def_hatphi}  the last substep of the fractional-step $\vartheta$-scheme may be written in the following compact form 
\begin{equation}
\label{compact_form}
\frac{U^n - \varPi^nU^{n-1}}{k_n} + \widehat{\varTheta}(t^{n-\frac 12}) =  P_0^n\hat\varphi (t^{n-\frac 12}).
\end{equation}

\section{Space-time reconstructions}

As aforementioned, the a posteriori error estimates will be derived by using the reconstruction technique. Our goal is  to define a continuous representation $\hat\omega$ of the approximate solution, $\hat{\omega}:[0,T]\rightarrow H_0^1(\varOmega),$ which will be a second order approximation of the exact solution $u(t)$ and whose residual will also be second order accurate. To define $\hat{\omega}(t)$ we shall exploit both the ideas of {\it elliptic reconstruction} introduced  in \cite{MN2003} and the fractional-step $\vartheta$-reconstruction based on approximations on one time subinterval introduced in  \cite{K2012}. Additionally, we shall extend the idea of the three-point Crank--Nicolson reconstruction introduced in  \cite{LPP2009, Prachittham2009}  and shall define a second fractional-step $\vartheta$-reconstruction  which will be based on approximations on  two adjacent  time  subintervals.  We shall begin our discussion by recalling the definition of the elliptic reconstruction operator and its basic properties.

\subsection{Reconstruction in space}
To derive a posteriori error estimates of optimal order in $L^\infty(0,T;L^2(\varOmega))$ norm for finite element discretizations of parabolic equations, the use of the 
elliptic reconstruction is necessary.  The elliptic reconstruction may be regarded as an a posteriori analogue to the Ritz--Wheeler projection appearing in standard a priori error
analysis for parabolic problems, c.f., for example, \cite{Wheeler73}, \cite{Thomee97}.  
Note that in fully discrete case with finite element spaces allowed to change with time 
the elliptic reconstruction operator depends on $n.$  
\begin{definition}[Elliptic reconstruction]  For fixed $v_n\in {\Bbb{V}}_h^n$, 
we define the elliptic 
reconstruction $\rec^n v_n\in H_0^1$ of $v_n,$ as the solution of the following elliptic problem
\begin{equation}
\label{ell_rec} 
\<\nabla \rec^n v_n, \nabla \psi\>=\<(-\varDelta_h^n) v_n,\psi\>\quad \forall \psi\in H_0^1.
\end{equation}
\end{definition}
 
It can be easily seen that the elliptic reconstruction $\rec^n$ satisfies the Galerkin orthogonality property
\begin{equation}
\label{ortho_pro}
\<\nabla(\rec^n v_n -v_n), \nabla \chi_n\>=0  \quad \forall \chi_n \in {\Bbb{V}}_h^n.
\end{equation}  

For completeness we shall  next give a residual-based a posteriori estimate for the elliptic reconstruction error $\|(\rec^n - I) v_n\|.$    \begin{lemma}[Residual-based a posteriori estimate for the elliptic reconstruction error] Let $v_n \in {\Bbb{V}}_h^n$ and $\rec^n v_n$ its elliptic reconstruction defined as in \eqref{ell_rec}. Then, it holds
\begin{equation}
\label{ell_rec_est}
\|(\rec^n - I) v_n\| \leq \eta_n(v_n),
\end{equation}
where $\eta_n$ is the elliptic estimator given by
\begin{equation}
\label{eta_n}
\begin{aligned}
\eta_{n}(v_n):=  C_{1,2}\, \|h_n^2(\varDelta-\varDelta_h^n)v_n\|_{\mathcal{T} _n}
+ C_{2,2} \, \|h_n^{3/2}J[\nabla v_n]\|_{\varSigma_n}\, .
\end{aligned}
\end{equation}
\end{lemma}
\proof 
Let $z \in H_0^1$ be the solution of problem 
\begin{equation}
\label{ell-dual}
a(\chi, z)=\<(\rec^n - I) v_n, \chi\> \quad \forall \chi\in  H_0^1, 
\end{equation}
and  $\mathcal{I}^nz\in\Bbb{V}_h^n$ be a Cl{\'e}ment-type interpolant of $z.$  
By using \eqref{ell-dual}, 
the orthogonality property  of the elliptic 
reconstruction \eqref{ortho_pro} and integration by parts, we get 

\begin{equation}
\begin{aligned}
\label{space_est_eq}
\|(\rec^n-I)v_{n}\|^2&=
a((\rec^n-I)v_{n} ,(z-\mathcal{I}^nz))\\  
&=\sum_{K\in T_n} \int_{K}(\varDelta-\varDelta_h^n)v_{n} (z-\mathcal{I}^nz)-
\sum_{e\in \varSigma_n } \int_{e}J[\nabla v_{n}] \,(z-\mathcal{I}^nz).
\end{aligned}
\end{equation}
Now, by applying first the approximation properties of the Cl\'ement interpolant \eqref{clement_app1} and afterwards the elliptic regularity \eqref{ell_reg}, we obtain
\begin{equation*}
\begin{aligned}
\sum_{K\in T_n} \int_{K}(\varDelta-\varDelta_h^n)v_{n} (z-\mathcal{I}^nz) 
&\le c_{1,2}\, |z|_2 \,\left(\sum _ {K \in \mathcal{T} _n} \|h_K^2(\varDelta-\varDelta_h^n)v_n\|_{K} ^2 \right)^{\frac 12} 
\\ & \le C_{1,2} \, \|(\rec^n-I)v_n\|\, \|h_n^2(\varDelta-\varDelta_h^n)v_n\|_{\mathcal{T} _n}
\end{aligned}
\end{equation*}
and
\begin{equation*}
\begin{aligned}
\sum_{e\in \varSigma_n } \int_{e}J[\nabla v_{n}] \,(z-\mathcal{I}^nz) 
&\le c_{2,2}\, |z|_2 \,\left(\sum _ {e\in \varSigma_n } \|h_K^\frac 32 \, J[\nabla v_{n}]\|_{e} ^2 \right)^{\frac 12}  
\\ & \le C_{2,2} \, \|(\rec^n-I)v_n\| \, \|h_n^{3/2}J[\nabla v_n]\|_{\varSigma_n}\, .
\end{aligned}
\end{equation*} \qed

We shall now turn our discussion to the time discretization and the so-called fractional-step $\vartheta$-reconstruction.

\subsection{Reconstruction in time}

Regarding the temporal variable, our goal is  to define a second order approximation $U(t)$ of $u(t),$ 
for all $t\in[0,T],$ whose residual is also second order accurate.  Choosing $U:[0,T]\rightarrow H_0^1(\varOmega)$ 
to be the piecewise linear interpolant at the nodal values, that is 
\begin{equation}\label{def_U}
U(t):=\ell_0^n(t)\, U^{n-1}+ \ell_1^n(t) \,U^{n}, \quad t\in I_n,
\end{equation}
where $\ell_0^n$ and $\ell_1^n$ are defined in \eqref{def_l0_l1}, 
 seems natural  for a second-order accurate scheme.  Indeed,  since the error at the nodes is of 
second order, $U(t)$ is an approximation of $u(t)$ of the same order, for all $t\in[0,T].$ 
However, its residual $R_U(t)$ 
\begin{equation}\label{res_U}
R_U(t):=U_t(t)-\varDelta U(t)-f(t), \quad t\in I_n,
\end{equation}
is  an a posteriori quantity of first order with respect 
to time.
We observe, using  \eqref{heat_eq}, that $R_U(t)$  may be written also in the form
\begin{equation}
\begin{aligned}
R_U(t)=[U_t(t)-u_t(t)]-\varDelta [U(t)-u(t)],\quad t\in I_n\, .
\end{aligned}
\end{equation}  
Although the second term on the right-hand side is 
of second order, we note that the first term is of first order only. 
By applying energy techniques to this 
error equation we can only derive residual-based a posteriori error estimates of suboptimal order with respect to time.

To recover the second order of accuracy in time,
we shall define appropriate reconstructions 
$\hat{U}$ in time which will be piecewise quadratic polynomials based on approximations on one 
time subinterval as well as on approximations based on two time subintervals.  
\begin{definition}[Time reconstructions] 
\label{time_rec}
We introduce
the piecewise quadratic time reconstruction  
$\hat{U}:[0,T]\rightarrow H_0^1(\varOmega),$  as follows 
\begin{equation}
\label{def_hat_U}
\hat{U}(t) = U(t) + \frac 12 (t-{t^{n-1}})(t-t^n) w_n, \quad t\in I_n,
\end{equation}
$n=1,\ldots, N,$ where $w_n$ is an appropriate piecewise constant polynomial with respect to time.
\end{definition}

In view of  \eqref{time_rec},  we can easily see that
\begin{lemma}[$L^\infty(L^2)$-estimate for the time reconstruction 
error] For $n = 1,\ldots, N,$ the following estimate holds
\begin{equation}
\label{time_rec_est}
\max_{t^{n-1}\leq t\leq t^n} \|(\hat{U}-U)(t)\| \leq  \frac{k_n^2}{8} \|w_n\|, \quad t\in I_n.
\end{equation}
\end{lemma}

In the sequel we shall study two  choices for the time reconstruction $\hat{U}$ which correspond to two appropriate choices for $w_n.$ In particular, we shall consider the following cases:
\begin{description}
\item[\emph{Time reconstruction 1}\rm{\emph{(based on approximations on one time subinterval)}}]
We shall extend the idea of the fractional-step $\th$-reconstruction introduced in  \cite{K2012} to the fully discrete case. For this purpose we choose $w_n$ as
\begin{equation}
 \label{def_wn}
w_n(t):= \varTheta_t(t) - P_0^n \varphi_t(t)= \frac{(-\varDelta_h^n) U^{n}-\widetilde\varPi^{n}(-\varDelta_h^{n-1}) U^{n-1}}{k_{n}}
-\frac{ P_0^n[f^{n}- f^{n-1}]}{k_{n}}, \quad t\in I_n.
\end{equation}
\item[\emph{Time reconstruction 2}\rm{\emph{(based on approximations on two adjacent  time subintervals)}}] 
The so-called three-point quadratic reconstruction for the Crank--Nicolson scheme,  \cite{LPP2009, Prachittham2009}, is defined  by choosing $w_n$ to 
be a finite difference approximation of $u_{tt}$ that uses the approximations on two time subintervals.  
Based on this idea,  we
define a three time-level quadratic reconstruction for the GFS-scheme by replacing $w_n$ in \eqref{def_hat_U} with 
\begin{equation}
\label{def_twn}
\widetilde{w}_n:= - {  \frac 2 {k_n + k_ {n-1} }
}
\Big [ \Bigl(\frac{U^{n}-\varPi^{n}U^{n-1}}{k_{n}}\Bigr)
-\pi^n\Bigl(\frac{U^{n-1}-\varPi^{n-1}U^{n-2}}{k_{n-1}}
\Bigr)\, \Big ], \quad t\in I_n,
\end{equation}
where $\pi ^n$ is any projection to $\Bbb{V}_h^n$ at our disposal.
\end{description}

\subsection{Reconstruction in both space and time}
The construction of appropriate space-time reconstructions for our analysis combines  the ideas discussed in the previous two paragraphs.  Let  $\omega:[0,T]\rightarrow H_0^1$ be the piecewise linear 
in time function  defined by linearly interpolating between the values $\rec^{n-1}U^{n-1}$ and $\rec^nU^n,$ 
\begin{equation}
\label{def_omega} 
\omega(t):=\ell_0^n(t)\rec^{n-1} U^{n-1}+\ell_1^n(t) \rec^nU^n, \quad t\in I_n,
\end{equation}
with $\ell_0^n$ and $\ell_1^n$ as in \eqref{def_l0_l1}.  According to the discussion above, the use of $\omega$ as intermediate 
function in order to derive  a posteriori error estimates will give optimal order estimates with respect the spatial derivative, however
it will lead to sub-optimal error estimates with respect the temporal one.  
The introduction of a piecewise quadratic polynomial in time is necessary, therefore we define the space-time reconstruction $\hat\omega$ of the approximate solution $U$: 
\begin{definition}[Space--time reconstruction] 
\label{space_time_rec}
  We introduce the space-time reconstruction $\hat\omega:[0,T]\rightarrow H_0^1$ of the approximate solution $U$  as follows
\begin{equation}
\label{def_hat_omega}
\hat \omega(t)= \omega(t) + \frac 12 (t-t^{n-1})(t^n-t)\rec^n w_n.
\end{equation}
\end{definition}

We shall now derive an $L^\infty(L^2)$-estimate for the space-time reconstruction error $\varepsilon = \hat \omega - U.$  The error $\varepsilon$ may be written as the sum of  the \emph{elliptic reconstruction error}  $\epsilon$  and  the \emph{time reconstruction error} $\sigma,$ that is    
\begin{equation}
\hat \omega- U = \epsilon + \sigma, \quad \text{where } \epsilon:=\omega -U, \;  \sigma:= \hat \omega -\omega.
\end{equation}

According to \eqref{ell_rec_est} and \eqref{time_rec_est} the following upper bounds for the reconstruction errors $\epsilon$ and $\sigma$ are valid.
\begin{lemma}For $m=1,\ldots,N,$ the following estimates hold
\begin{equation}
\max_{0\leq t\leq t^m} \|\epsilon(t)\| \leq \E_m^{\emph{ell}} \quad \text{with} \quad \E_m^{\emph{ell}}:=\max_{0\leq n \leq m} \eta_{n}(U^{n})
\end{equation}
and
\begin{equation}
\max_{0\leq t\leq t^m} \|\sigma(t)\| \leq   \E_m^{\emph{rec}}(w_n)  \quad \text{with} \quad \E_m^{\emph{rec}}(w_n):=\max_{0\leq n \leq m}\frac{k_n^2}{8} \bigl\{\eta_n(w_n)
+ \|w_n\|\bigr\},
\end{equation}
where $\eta_n$ is defined in \eqref{eta_n}.
\end{lemma}

\section{$L^\infty(L^2)$-estimates for the total error}
Let   $\rho,\; \hat\rho$  denote  the {\it parabolic errors} defined by $\rho:=u-\omega,\;\hat\rho:=u-\hat\omega,$ respectively.
The \emph{total error} $e:=u-U$ can be  split as follows 
\begin{equation}
e=u-U = (u- \hat \omega) +(\hat \omega - U) =
\hat\rho + \sigma + \epsilon.
\end{equation}
A bound for the reconstruction error $\sigma +\epsilon$ was presented in the previous section.  We shall now continue with the estimation of  the basic parabolic error,    which is  
stated in Theorem \ref{Theorem3_1}.

\subsection{An a posteriori estimate in $L^\infty(L^2)$ and $L^2(H^1)$ for 
the parabolic error}

We  begin with the derivation of  the error equation: 
\begin{lemma}[Error equation]
\label{Lemma2.1}
For each $t\in I_n,$ it holds
\begin{equation}
\label{eq_err1}
\<\hat\rho_t(t),\psi\>+a(\rho(t),\psi)= \<R_h, \psi\>  \quad \forall \psi\in H_0^1,
\end{equation}
with
\begin{equation}
\label{def_Rh}
\begin{aligned}
R_h:=& \ell_0^n(t) (\varPi^n-I)(-\varDelta_h^{n-1})U^{n-1} -k_n^{-1}(\varPi^n -I)U^{n-1} -(t-t^{n-\frac 12}) (\rec^n-I) w_n\\&-(t-t^{n-\frac 12}) (w_n -\varTheta_t(t) +P_0^n\varphi_t (t))
 - \frac{(\rec^{n-1}-I) U^{n-1}-(\rec^{n}-I)U^{n}}{k_n} \\&
 + \xi_\varTheta^n + f(t)- P_0^n\hat\varphi(t).
\end{aligned}
\end{equation}
\end{lemma}
\proof 
For $n=1,\ldots,N,$ and  $\psi\in H_0^1,$ in view of 
\eqref{def_omega} and \eqref{def_hat_omega}, we obtain
\begin{equation}
\label{err_theta1}
\begin{aligned}
\<\hat \omega_t(t),\psi\>= \<\omega_t(t), \psi\> + (t-t^{n-\frac 12}) \<\rec^n w_n, \psi\>.
\end{aligned}
\end{equation}

Thus, by using Definition \ref{space_time_rec}, we get
\begin{equation*}
\begin{aligned}
\<\hat\rho_t(t),\psi\>+a(\rho(t),\psi)=
&-a(\ell_0^n(t) \rec^{n-1}U^{n-1}+\ell_1^n(t)\rec^n U^n, \psi)
-(t-t^{n-\frac 12}) \<\rec^n w_n, \psi\>\\&-
k_n^{-1}\<\rec^n U^{n}-\rec^{n-1}U^{n-1},\psi\>+
\<f(t),\psi\> %
\end{aligned}
\end{equation*}
According to the elliptic reconstruction definition \eqref{ell_rec}, the last relation leads to
\begin{equation}
\label{err_theta3}
\begin{aligned}
\<\hat\rho_t(t),&\psi\>+a(\rho(t),\psi)=-\< \ell_0^n(t)(-\varDelta_h^{n-1})U^{n-1} +\ell_1^n(t) (-\varDelta_h^{n})U^n,\psi\> 
\\&-(t-t^{n-\frac 12}) \<\rec^n w_n, \psi\> -
k_n^{-1}\<\rec^nU^{n}-\rec^{n-1}U^{n-1},\psi\>+ \<f(t),\psi\>,
\end{aligned}
\end{equation}
from which, in view of \eqref{def_Theta}, we infer that
\begin{equation}
\label{err_theta4}
\begin{aligned}
\<\hat\rho_t(t),\psi\>+a(\rho(t),\psi)= &-\<\varTheta(t), \psi\> + \ell_0^n(t)\<(\varPi^n-I)(-\varDelta_h^{n-1})U^{n-1}, \psi\> \\&-(t-t^{n-\frac 12}) \<(\rec^n-I) w_n , \psi\>-(t-t^{n-\frac 12}) \<w_n , \psi\>\\&
 -k_n^{-1}\<(\rec^{n-1}-I) U^{n-1}-(\rec^{n}-I)U^{n},\psi\> \\&+k_n^{-1}\<U^{n}-U^{n-1},\psi\>+ \<f(t),\psi\>,
\end{aligned}
\end{equation}
Now, in view of  \eqref{compact_form}, we observe that
\begin{equation}
\label{err_theta4}
\begin{aligned}
\<\hat\rho_t(t),&\psi\>+a(\rho(t),\psi)=  \ell_0^n(t)\<(\varPi^n-I)(-\varDelta_h^{n-1})U^{n-1} -k_n^{-1}(\varPi^n -I)U^{n-1},\psi\> \\&-(t-t^{n-\frac 12}) \<(\rec^n-I) w_n , \psi\>-(t-t^{n-\frac 12}) \<w_n , \psi\>+\<\hat\varTheta(t^{n-\frac 12})- P_0^n\hat\varphi (t^{n-\frac 12}),\psi\>\\&
 -k_n^{-1}\<(\rec^{n-1}-I) U^{n-1}-(\rec^{n}-I)U^{n},\psi\> -\<\varTheta(t), \psi\>+ \<f(t) ,\psi\> .
\end{aligned}
\end{equation}
In view of \eqref{def_hatTheta} and \eqref{def_hatphi} we can easily see that
\begin{equation}
\hat\varTheta(t^{n-\frac 12})- P_0^n\hat\varphi (t^{n-\frac 12}) = (t-t^{n-\frac 12})(\varTheta_t(t) - P_0^n\varphi_t (t) ),
\end{equation}
and thus to conclude that
\begin{equation}
\label{err_theta4}
\begin{aligned}
\<\hat\rho_t(t),\psi\>+&a(\rho(t),\psi)=  \ell_0^n(t)\<(\varPi^n-I)(-\varDelta_h^{n-1})U^{n-1} -k_n^{-1}(\varPi^n -I)U^{n-1},\psi\> \\&-(t-t^{n-\frac 12}) \<(\rec^n-I) w_n , \psi\>-(t-t^{n-\frac 12}) \<w_n - \varTheta_t(t) +P_0^n\varphi_t (t),\psi\>\\&
 -k_n^{-1}\<(\rec^{n-1}-I) U^{n-1}-(\rec^{n}-I)U^{n},\psi\> \\&+\<\xi_\varTheta^n, \psi\>+ \<f(t)- P_0^n\hat\varphi(t),\psi\> .
\end{aligned}
\end{equation}
\qed

An a posteriori error bound for the parabolic error follows. 
Note that the estimate that we will derive depends still on the choice of the time reconstruction 
through  $w_n$ as well as on stationary finite element errors through 
the elliptic reconstruction $\rec^n .$ 
  
\begin{theorem} {\sc (Estimates in $L^\infty(L^2)$ and $L^2(H^1)$ for 
the parabolic error)}
\label{Theorem3_1}

The following estimate is valid
\begin{equation}
\max_{t\in[0,t^m]}\bigl \{\|\hat\rho(t)\|^2 + 
  \int _0 ^{t} (|\rho(s)|^2 + |\hat\rho(s)|^2) \, ds\bigr\}  
\leq \,  \|\hat\rho(0)\|^2  + \mathcal{J}_m,
\end{equation}
where $\mathcal{J}_m, \, m=1,\ldots,N,$ is defined by 
\begin{equation}
\mathcal{J}_m:= \mathcal{J}_m^{T,1}
+ \mathcal{J}_m^{T,2}+ \mathcal{J}_m^{S,1} + \mathcal{J}_m^{S,2}+
 \mathcal{J}_m^{C}+ \mathcal{J}_m^{D} + \mathcal{J}_m^{W},
\end{equation}
with
\begin{equation}
\mathcal{J}_m^{T,1} := \sum_{n=1}^{m}\int_{t^{n-1}}^{t^n} |\sigma(s)|^2\,  ds,
\end{equation}
\begin{equation}
\mathcal{J}_m^{T,2} :=2 \sum_{n=1}^{m}\int_{t^{n-1}}^{t^n}| \<\xi_\varTheta^n,\hat\rho(s)\>|  ds,
\end{equation}
\begin{equation}
\mathcal{J}_m^{S,1} := 2\sum_{n=1}^{m}\int_{t^{n-1}}^{t^n}  
|(s-t^{n-\frac 12})\<(\rec^n-I)w_n, \hat\rho(s)\>| ds  
\end{equation}
\begin{equation}
\mathcal{J}_m^{S,2}:=2\sum_{n=1}^{m} \int _{t^{n-1}}^{t^n} |\<\frac{(\rec^n - I) U^{n}-(\rec^{n-1}-I )U^{n-1}}{k_n} ,\hat\rho(s)\>|ds \end{equation}
\begin{equation}
\mathcal{J}_m^{C}:=2\sum_{n=1}^{m}\int _{t^{n-1}}^{t^n}  
|\<\ell_0^n(s) (\widetilde\varPi^n-I)(-\varDelta_h^{n-1}) U^{n-1} -k_n^{-1}(\varPi^n-I) U^{n-1},\hat\rho(s) \>|\,ds\, .
\end{equation}
\begin{equation}
\mathcal{J}_m^{D}:=2\sum_{n=1}^{m}\int _{t^{n-1}}^{t^n}  
| \<f(s)- P_0^n\hat\varphi(s), \hat\rho(s)\>|\,ds
\end{equation}\begin{equation}
\mathcal{J}_m^{W}:=2\sum_{n=1}^{m}\int _{t^{n-1}}^{t^n}  
|(s-t^{n-\frac 12}) \<w_n - \varTheta_t(s) +P_0^n\varphi_t (s), \hat\rho(s)\>|\,ds
\end{equation}
\end{theorem}
\proof
Setting $\psi=\hat\rho$ in \eqref{eq_err1} and observing that
\begin{equation*}
a(\rho(t),\hat\rho(t))=\frac{1}{2}|\rho(t)|^2 + \frac{1}{2}|\hat\rho(t)|^2
-\frac{1}{2} |\hat\rho(t)-\rho(t)|^2,
\end{equation*}
we obtain
\begin{equation*}
\begin{aligned}
\|\hat\rho(t)\|^2  + 
  \int _0 ^{t} (|\rho(s)|^2 + |\hat\rho(s)|^2) \, ds 
\leq \|\hat\rho(0)\|^2 +  \int_0^{t^m} |\sigma(s)|^2\,  ds  +  2 \int_0^{t^m} \<R_h,\hat\rho(s)\>\, ds, 
\end{aligned}
\end{equation*}
for all $t \in [0, t^m].$ Recalling the 
definition \eqref{def_Rh} of $R_h,$  it can be easily seen that
\begin{equation}
2 \int_0^{t^m} \<R_h,\hat\rho(s)\>\, ds \leq  \mathcal{J}_m^{T,2}+\mathcal{J}_m^{S,1} +\mathcal{J}_m^{S,2}+
\mathcal{J}_m^{C}+\mathcal{J}_m^{D}+\mathcal{J}_m^{W},
\end{equation}
which completes the proof.
\qed

We emphasize here that the piecewise polynomial in time $w_n$ appearing in the definition of time reconstruction \eqref{def_hat_U}
is chosen such that $\mathcal{J}_m^{W}$ is an a  posteriori quantity of optimal order.  According to \eqref{def_wn},  the term $\mathcal{J}_m^{W}$ vanishes in case of the time reconstruction based on one 
time subinterval. In addition, in case of the time reconstruction based on two adjacent time subintervals the following result is valid:
\begin{lemma}[Calculation of $\widetilde{w}_n - \varTheta_t(t) +P_0^n\varphi_t (t)$]
\label{Lemma2.2} For $t\in I_n$ we have
\begin{equation}
\label{term_twn}
\widetilde{w}_n- \varTheta_t(t)+P_0^n\varphi_t(t) = - \frac {2} {k_n + k_ {n-1} }\bigl\{z_n +\xi_\varTheta^n - \pi^n  \xi_\varTheta^{n-1} - y_n -P_0^n\xi_\varphi^n  + \pi^n  P_0^{n-1}\xi_\varphi^{n-1}\bigr\}
\end{equation}
with $z_n$ and $y_n$ defined by
\begin{equation}
\label{def_zn_yn}
\begin{aligned}
z_n := & \frac{1} {2 }\, \Bigl [\ 
\frac{k _{n-1} }{   k_n}  
(-\varDelta_h^n)U^n\ - \big ( (2 + \frac{k _{n-1} }{  k_n}) \widetilde\varPi ^n - 
 \, \pi^n  \big ) \,
(-\varDelta_h^{n-1})U^{n-1} \\& \quad
 +   \,    \, \pi^n \, \widetilde\varPi^{n-1} (-\varDelta_h^{n-2})U^{n-2} 
\Bigr ],\\
y_n :=  &\frac{1} {2 }\,  \Bigl [\ 
\frac{k _{n-1} }{   k_n}  
P_0^nf^n\ - \big ((2  + \frac{k _{n-1} }{  k_n})P_0^n-\pi^n P_0^{n-1}\big )f^{n-1} 
 +  \pi^n P_0^{n-1} f^{n-2} 
\Bigr ].
\end{aligned}
\end{equation}\end{lemma}
\proof 
We let $\widetilde\varphi$ be given by
\begin{equation}
\label{def_tildephi}
\widetilde\varphi(t):= \ell_{1/2}^n(t) P_0^n\varphi(t^{n-\frac 12}) + \ell_{-1/2}^n(t)  \pi^n  P_0^{n-1}\varphi(t^{n-\frac 32}), \quad  t\in I_n, 
\end{equation}
where
\begin{equation}
\label{l12}
\ell_{1/2}^n(t):=\frac{2\,(t-t^{n-\frac 3 2})}{k_n+k_{n-1}}, 
\quad \ell_{-1/2}^n(t):=\frac{2\,(t^{n-\frac 1 2}-t)}{k_n+k_{n-1}} . 
\end{equation}
We express $\varTheta(t), \, \varphi(t),\, t\in I_n,$ defined in \eqref{def_Theta} and in \eqref{def_phi}, respectively, in terms of $\ell_{1/2}^n $ and $\ell_{-1/2}^n,$ that is
\begin{equation}
\begin{aligned}
\varTheta (t)=\,&\ell_{1/2}^n(t)\, \varTheta(t^{n-\frac 1 2}) 
+\ell_{-1/2}^n(t)\, \widetilde \varTheta^{n-\frac 3 2} , \quad t\in I_n,\\
\varphi (t)=\,&\ell_{1/2}^n(t)\, \varphi(t^{n-\frac 1 2}) 
+\ell_{-1/2}^n(t)\,\widetilde \varphi^{n-\frac 3 2}, \quad t\in I_n,
\end{aligned}
\end{equation}
where
\begin{equation}
\label{def_tildetheta}
\begin{aligned}
\widetilde \varTheta ^{n-\frac 3 2} :=\,
& \ell_{0}^n(t^{n-\frac 3 2} )\,  
 \widetilde\varPi^n (-\varDelta_h^{n-1})U^{n-1}
+\ell_{1}^n(t^{n-\frac 3 2})\,  (-\varDelta_h^n)U^n  \, ,\\
\widetilde \varphi ^{n-\frac 3 2} :=\,
& \ell_{0}^n(t^{n-\frac 3 2} )\,  
 \varphi(t^{n-1})
+\ell_{1}^n(t^{n-\frac 3 2})\,   \varphi(t^{n})\, .\end{aligned}
\end{equation}
Now, in view of \eqref{def_twn} and \eqref{compact_form}, we have that  
\begin{equation*}
\begin{aligned}
\widetilde{w}_n-\varTheta_t(t)+P_0^n\varphi_t(t) =&-\frac {2} {k_n + k_ {n-1} }\Bigl\{\Bigl(\frac{U^{n}-\widetilde\varPi^{n}U^{n-1}}{k_{n}}\Bigr)
 -\pi^n\Bigl(\frac{U^{n-1}-\widetilde\varPi^{n-1}U^{n-2}}{k_{n-1}}
\Bigr)\Bigr\}
\\&\quad-\varTheta _t(t) +P_0^n\varphi_t(t)\\
 = &-\frac {2} {k_n + k_ {n-1} }\Bigl\{ -\widehat\varTheta(t^{n-\frac 12}) + P_0^n\hat\varphi(t^{n-\frac 12})
- \pi^n[- \widehat\varTheta(t^{n-\frac 32}) + P_0^n\hat\varphi(t^{n-\frac 32})]\\& \quad+ \varTheta(t^{n-\frac 1 2}) 
-\widetilde\varTheta^{n-\frac 3 2} - P_0^n \varphi(t^{n-\frac 12}) + P_0^n \widetilde\varphi^{n-\frac 32} \Bigr\}\, .
\end{aligned}
\end{equation*}
According to \eqref{def_hatTheta} and \eqref{def_hatphi}, we get
\begin{equation}
\label{data_term_rec2}
\begin{aligned}
\widetilde{w}_n-\varTheta_t(t)+P_0^n\varphi_t(t) &=   - \frac {2} {k_n + k_ {n-1} }\Bigl\{ 
\pi^n \varTheta(t^{n-\frac 32}) -\widetilde\varTheta^{n-\frac 3 2} +\xi_\varTheta^n - \pi^n  \xi_\varTheta^{n-1}
\\&-\pi^n P_0^{n-1}\varphi(t^{n-\frac 32}) + P_0^n \widetilde\varphi^{n-\frac 32} -P_0^n\xi_\varphi^n  + \pi^n  P_0^{n-1}\xi_\varphi^{n-1}\Bigr\}, \quad  t\in I_n
.\end{aligned}
\end{equation}
In view of \eqref{def_Theta}, \eqref{def_tildetheta} and \eqref{def_zn_yn}, we can easily see that
\begin{equation}
\pi^n \varTheta(t^{n-\frac 32}) -\widetilde\varTheta^{n-\frac 3 2} = z_n \quad \text{and} \quad  
\pi^n P_0^n\varphi(t^{n-\frac 32}) - P_0^n \widetilde\varphi^{n-\frac 32} = y_n,
\end{equation}
and the desired result follows.
\qed 

Note that, in case of constant time-steps and mesh,
$z_n$ corresponds to a $k^2(-\varDelta) u_{tt}$ term of optimal order.

In the next section, we shall further  investigate 
each term of the estimator  by considering both time reconstructions 
in combination with residual-based a posteriori estimators for the elliptic error; 
other choices of estimators for the stationary finite element errors are also possible.  %
\smallskip


\subsection{A residual-based a posteriori bound for the parabolic error.}
\label{sec:errest}

In this paragraph we use the space-time reconstruction introduced in \eqref{def_hat_omega}, with $w_n$ to be chosen
either as in \eqref{def_wn} or as in \eqref{def_twn}, and 
residual-based estimators to derive an upper bound for the parabolic error $\hat\rho.$    The proof is split in several steps.

Throughout the rest of this paragraph we denote by  $t^m_\star \in [0,t^m]$ the time for which
\begin{equation}
\|\hat\rho(t_\star^m)\|= \max_{t\in [0, t^m]} \|\hat\rho(t)\|.
\end{equation}

We shall first show an upper  bound for the terms $\mathcal{J}_m^{T,1}$ and $\mathcal{J}_m^{T,2}$ 
appearing in Theorem \ref{Theorem3_1},  
which measure the local time discretization error.  

\begin{lemma}[Time error estimate]
Let  $v_n \in {\Bbb{V}}_h^n$ and the time estimators $\E_m^{T,1}$  and $\E_m^{T,2}$ be defined as follows
\begin{equation}
\E^{T,1}_m(v_n): = \left(\sum_{n=1}^{m} k_n \,\gamma_n^2(v_n)\right)^{1/2} \quad \text{with} \quad \gamma_n(v_n) := \frac{k_n^2}{\sqrt{30}} \bigl (c_1 |v_n| +  C_{1,1}\|h_n(-\varDelta_h^n)v_n\| \bigr)\,  
\end{equation}
and
\begin{equation}
\E^{T,2}_m:= 2\sum_{n=1}^{m} k_n\,   \|\xi_{\Theta}^n\|.
 \end{equation}
Then, we have
\begin{equation}
\begin{aligned}
\mathcal{J}_m^{T,1} \leq \left(\E^{T,1}_m(w_n)\right)^2 \quad \text{or} \quad \mathcal{J}_m^{T,1} \leq  \left(\E^{T,1}_m(\widetilde{w}_n)\right)^2,
\end{aligned}
\end{equation}
depending on the choice of the time reconstruction  with $w_n$ and $\widetilde{w}_n$ be defined  in \eqref{def_wn} and \eqref{def_twn}, respectively. In addition,  the term $\mathcal{J}_m^{T,2}$  may be bounded as follows
\begin{equation}
\mathcal{J}_m^{T,2}\leq  \|\hat \rho(t^m_{\star})\|\, \E^{T,2}_m.
\end{equation}
\end{lemma}
\proof
We have
\begin{equation}
\label{time_est_rel}
|\sigma(t)|_1^2 = a(\sigma(t),\sigma(t))
= \frac 12 (t-t^{n-1})(t^n-t)a(\rec^nw_n, \sigma(t)).
\end{equation}
By using the definition of the elliptic reconstruction \eqref{ell_rec}, we get
\begin{equation*}
\begin{aligned}
|\sigma(t)|_1
&\leq (t-t^{n-1})(t^n-t)|(-\varDelta_h^n)w_n|_{-1}.
\end{aligned}
\end{equation*}
To estimate the dual norm in the above relation,  we can proceed as follows
\begin{equation}
\label{time_est_rel2}
\begin{aligned}
|-\varDelta_h^n w_n|_{-1}&=\sup_{0\neq z\in H_0^1} 
\frac{\<-\varDelta_h^n w_n, z\>}{|z|_1}\\&=
\sup_{0\neq z\in H_0^1}\left \{ 
\frac{\<-\varDelta_h^nw_n, \mathcal{I}^nz\>}{|z|_1} 
+ \frac{\<-\varDelta_h^nw_n, z-\mathcal{I}^nz\>}{|z|_1}\right \},
\end{aligned}
\end{equation}
with  $\mathcal{I}^nz\in \Bbb{V}_h^n$ a Cl{\'e}ment-type interpolant of  $z.$
Now, in view of \eqref{dis_laplacian} and \eqref{clement_stab}, we have
\begin{equation}
\label{time_est_rel3}
\<-\varDelta_h^n w_n,\mathcal{I}^nz\> 
\leq c_1|w_n|_1 |z|_1. 
\end{equation}
Furthermore, using the approximation properties  \eqref{clement_app1} of a Cl{\'e}ment-type 
interpolant, we obtain
\begin{equation}
\label{time_est_rel4}
\<-\varDelta_h^nw_n, z- \mathcal{I}^nz\> \leq
c_{1,1}\|h_n (-\varDelta_h^n) w_n\|\, |z|_1.
\end{equation} 
According to \eqref{time_est_rel3} and \eqref{time_est_rel4}, \eqref{time_est_rel2} leads to  
\begin{equation}
\begin{aligned}
|-\varDelta _h^n w_n|_{-1} \leq  
c_1  |w_n|_1 
+ C_{1,1}\|h_n(-\varDelta_h^n) w_n\|.
\end{aligned}
\end{equation}
By observing that
\begin{equation}
\label{time_est_rel1}
\begin{aligned}
\int_{t^{n-1}}^{t^n}(t-t^{n-1})^2(t^n-t)^2 dt  =\frac{k_n^5}{30} 
\end{aligned}
\end{equation}
the desired result follows.
\qed

We shall next estimate the term $\mathcal{J}_m^{S,1}$ in Theorem \ref{Theorem3_1} which accounts
for the space discretization error. 

\begin{lemma}[Spatial error estimate] Let  $v_n \in {\Bbb{V}}_h^n$ and $\E_m^{S,1}$ be defined as
\begin{equation*}
 \E^{S,1}_m(v_n):= \sum_{n=1}^{m} \frac{k_n^2}{2} \eta_n(v_n).
\end{equation*}

Then, depending on the choice of the time reconstruction,  the following estimate is valid
\begin{equation}
 \mathcal{J}_m^{S,1}\leq  \|\hat \rho(t^m_{\star})\| \, \E^{S,1}_m(w_n)
\quad \text{or} \quad  \mathcal{J}_m^{S,1}\leq  \|\hat \rho(t^m_{\star})\| \,\E^{S,1}_m( \widetilde{w}_n)\, ,
\end{equation}
where $w_n$ and $\widetilde{w}_n$ are defined in \eqref{def_wn} and \eqref{def_twn}, respectively.\end{lemma}
\proof
Since $w_n$  is piecewise constant in time, we can easily see that
\begin{equation}
\int_{t^{n-1}}^{t^n}  
|(s-t^{n-\frac 12})\<(\rec^n-I)w_n, \hat\rho(s)\>| ds  \leq  \,  \frac {k_n^2} 4\,  \max_{t\in [t^{n-1}, t^n]} \|\hat\rho(t)\| \,  \|(\rec^n-I)w_n\|
\end{equation}
and the desired result follows.
\qed
\smallskip

An upper bound for  the term $\mathcal{J}_m^{S,2}$ in Theorem \ref{Theorem3_1} will be next presented. 

\begin{lemma}[Space estimator accounting for mesh changing] Let $\E_m^{S,2}$ be defined as
\begin{equation}
\E^{S,2}_m:=2\sum_{n=1}^{m} k_n   \delta_n 
\end{equation}
with
\begin{equation}
\begin{aligned}
\delta_n:=&\bigl\{
C_{1,2}\|\check{h}_n^2\big  [k_n^{-1}(\varDelta-\varDelta_h^n)U^n- 
k_n^{-1}(\varDelta-\varDelta_h^{n-1})U^{n-1}\big ]\|_{\check{\mathcal{T}} _n}
 \\&+
C_{2,2} \|\check{h}_n^{3/2} 
J[\nabla U^{n}-\nabla U^{n-1}]\|_{\check{\varSigma}_n} \bigr\}.  
\end{aligned}
\end{equation}
Then, we have that
\begin{equation}
\mathcal{J}_m^{S,2}\leq \|\hat \rho(t^m_{\star})\| \,   \E_m^{S,2}\, .
\end{equation}
\end{lemma}
\proof
Let $z:[0,T]\rightarrow H_0^1$ be the solution of problem 
\begin{equation}
\label{eq:dual}
a(\chi, z(t))=\<\hat\rho(t), \chi\>, \quad \forall \chi\in  H_0^1, \; t\in[0,T],
\end{equation}
and  $\hat{\mathcal{I}}^nz(t)\in\Bbb{V}_h^n\cap \Bbb{V}_h^{n-1},\; t\in I_n,$ be its Cl{\'e}ment-type interpolant. 
Since  $\hat{\mathcal{I}}^nz(t)\in\Bbb{V}_h^n\cap \Bbb{V}_h^{n-1},$ using first \eqref{eq:dual}, 
the orthogonality property  of the elliptic 
reconstruction \eqref{ortho_pro} in $\Bbb{V}_h^{n-1}\cap\Bbb{V}_h^{n},$ and integration by parts, we get 
\begin{equation}
\begin{aligned}
\label{space_est_eq}
\<(\rec^n-I)U^{n}-&(\rec^{n-1}-I )U^{n-1},\hat\rho(t)\>\\&=
a((\rec^n-I)U^{n}-(\rec^{n-1}-I )U^{n-1} ,(z-\hat{\mathcal{I}}^nz)({t}))\\  
&=\sum_{K\in \check{T}_n} \int_{K}\bigl \{(\varDelta-\varDelta_h^n)U^{n} - 
(\varDelta - \varDelta_h^{n-1})U^{n-1}\bigr\} (z-\hat{\mathcal{I}}^nz)(t)\\&\qquad-
\sum_{e\in \check{\varSigma}_n } \int_{e}J[\nabla U^{n}-\nabla U^{n-1}] \,(z-\hat{\mathcal{I}}^nz)({t})  \bigl\}.
\end{aligned}
\end{equation}

Hence, in view of \eqref{clement_app1}, we obtain
\begin{equation}
\begin{aligned}
\sum_{K\in \check{T}_n} \int_{K}\bigl \{(\varDelta-\varDelta_h^n)U^{n}&- (\varDelta - \varDelta_h^{n-1})U^{n-1}\bigr\} (z-\hat{\mathcal{I}}^nz)(t) \\ &\leq 
C_{1,2} \|\check{h}_n^2\bigl \{(\varDelta-\varDelta_h^n)U^n-(\varDelta-\varDelta_h^{n-1})U^{n-1}\bigr\}\|\, \|\hat\rho(t)\|,
\end{aligned}
\end{equation} 

and
\begin{equation}
\sum_{e\in \check{\varSigma}_n } \int_{e}J[\nabla U^{n}-\nabla U^{n-1}] 
(z-\hat{\mathcal{I}}^nz)(t) 
 \leq  C_{2,2} \|\check{h}_n^{3/2} 
J[\nabla U^{n} - \nabla U^{n-1}]\|_{\check{\varSigma}_n}\, \|\hat\rho(t)\|;
\end{equation} 
the claimed result follows.
\qed

The term $\mathcal{J}_m^{C}$ in Theorem \ref{Theorem3_1} that corresponds to the coarsening error can be bounded as follows
\begin{lemma}[Coarsening error estimate] Let $\E_m^C$ be the coarsening estimator defined by
\begin{equation}
\E_m^C := 2 \sum_{n=1}^{m}k_n\beta_n \quad \text{with} \qquad \beta_n:= \|(\varPi^n-I)(\varDelta_h^{n-1} U^{n-1} +k_n^{-1}U^{n-1} )\|_{\mathcal{T} _n}.
\end{equation}
Then, it holds
\begin{equation}
\mathcal{J}_m^{C}\leq \|\hat \rho(t^m_{\star})\|\, \E_m^C.
\end{equation}
\end{lemma}

Upper bounds for the term $\mathcal{J}_m^{D}$  which measure the data approximation error,
will be shown in the next lemma. 
\begin{lemma}[Data error estimate] Let 
\begin{equation}
\label{zeta1_n}
\begin{aligned}
\zeta_{n,1}&:= \frac{1}{k_n} \int_{t^{n-1}}^{t^n} \|f(s)-\varphi(s)\| \, ds, \\
\zeta_{n,2}&:= c_{1,1}\max \bigl\{\|h_n(I-P_0^n)(f^{n-1} + \xi_\varphi^n)\|, \|h_n(I-P_0^n)(f^{n}+\xi_\varphi^n)\|\bigr\},
\end{aligned}
\end{equation}
and
\begin{equation}
\E_m^{D,1} := 2\, \sum_{n=1}^{m} k_n \, (\zeta_{n,1} + \|\xi_\varphi^n\|)
 \end{equation}
with $\xi_\varphi^n$ defined in  \eqref{corr_phi}.
Then,  we have that
\begin{equation}
\begin{aligned}
\mathcal{J}_m^{D,2} &\leq  \|\hat \rho(t^m_{\star})\|\,\E_m^{D,1} +  \,2\sum_{n=1}^{m}  \left(\int_{t^{n-1}}^{t^n}|\hat\rho(s) |^2\right)^{1/2}k_n^{1/2} \zeta_{n,2}\, .
\end{aligned}
\end{equation}
\end{lemma}
\proof The term $\mathcal{J}_m^{D,2}$ may be  
bounded as follows
\begin{equation}
\begin{aligned}
\mathcal{J}_m^{D,2} = &2\,\sum_{n=1}^{m} \int_{t^{n-1}}^{t^n} |\<f(s)-P_0^n\hat\varphi(s), \hat\rho(s)\>|\,ds \\
\leq &2\,\sum_{n=1}^{m} \int_{t^{n-1}}^{t^n} \bigl \{|\<f(s)-\hat\varphi(s), \hat\rho(s)\>| +
|\<(I-P_0^n)\hat\varphi(s), \hat\rho(s)\>|\bigr\} \, ds .
\end{aligned}
\end{equation}
Now, we have
\begin{equation*}
\begin{aligned}
&\int_{t^{n-1}}^{t^n}  
|\<f(s)-\hat\varphi(s), \hat\rho(s)\>| \, ds \leq  \max_{t\in[t^{n-1}, t^n]} \|\rho(t)\| \,\int_{t^{n-1}}^{t^n} \bigl \{\|f(s)-\varphi(s)\|+  \|\xi_\varphi^n\| \bigr\} \, ds,
\end{aligned}
\end{equation*}
from which we can conclude that
\begin{equation}
\label{JD2_1}
\begin{aligned}
&\sum_{n=1}^{m} \int_{t^{n-1}}^{t^n}  
|\<f(s)-\hat\varphi(s), \hat\rho(s)\>| \, ds \leq \|\hat \rho(t^m_{\star})\|\,\sum_{n=1}^{m}  k_n \, (\zeta_{n,1} + \|\xi_\varphi^n\|) \, .
\end{aligned}
\end{equation}
Furthermore,  using again the orthogonality property of $P_0^n$, we obtain 
\begin{equation*}
\begin{aligned}
\<(I-P_0^n)\hat\varphi(s),\hat\rho(s) \> = \<(I-P_0^n)\hat\varphi(s),
(\hat\rho-\mathcal{I}^n \hat\rho )(s)\> 
\leq c_{1,1}\|h_n (I-P_0^n)\hat\varphi(s)\| \, |\hat\rho(s) |_1 .
\end{aligned}
\end{equation*}
Now,
\begin{equation*}
\begin{aligned}
\|h_n (I-P_0^n)\hat\varphi(s)\|& = \| h_n(I-P_0^n)[\ell_0^n(s) f^{n-1}
+ \ell_1^n(s) f^{n} + \xi_\varphi^n]\|\\
&\leq\max \bigl\{\|h_n(I-P_0^n)(f^{n-1} + \xi_\varphi^n)\|, \|h_n(I-P_0^n)(f^{n}+\xi_\varphi^n)\|\bigr\} ,
\end{aligned}
\end{equation*}
and hence, 
\begin{equation}
\label{JD2_2}
\begin{aligned}
&\sum_{n=1}^{m}  \int_{t^{n-1}}^{t^n} |\<(I-P_0^n)\hat\varphi(s),\hat\rho(s) \>| \, ds 
\leq \sum_{n=1}^{m}  c_{1,1}\,\left(\int_{t^{n-1}}^{t^n}|\hat\rho(s) |^2\right)^{1/2}  k_n^{1/2} \,
  \zeta_{n,2}  .
\end{aligned}
\end{equation}
In view of \eqref{JD2_1}, \eqref{JD2_2}, we conclude the desired result.
\qed

\begin{lemma}[An estimator for $\mathcal{J}_m^{W}$] The term $\mathcal{J}_m^{W}$ vanishes in case of the two time-level reconstruction.  Furthermore,  the term $\mathcal{J}_m^{W}$ corresponding to the three time-level reconstruction may be bounded as follows
\begin{equation}
\begin{aligned}
\mathcal{J}_m^{W} &\leq \, \|\hat \rho(t^m_{\star})\| \,\sum_{n=1}^{m} \frac{k_n^2}{2(k_n +k_{n-1})} \,  \|z_n\| + \|\hat \rho(t^m_{\star})\| \,\sum_{n=1}^{m}\frac{ k_n^2}{4}  \left(\|\xi_\Theta^n\| + \|\pi^n  \xi_\varTheta^{n-1}\| \right)
\\&+\|\hat \rho(t^m_{\star})\|\,\sum_{n=1}^{m} \frac{k_n^2}{2(k_n +k_{n-1})}  \|y_n\| +\|\hat \rho(t^m_{\star})\| \,\sum_{n=1}^{m}\frac{ k_n^2}{4} \,  \left( \|P_0^n\xi_\varphi^n\|
 +\ \|\pi^n P_0^{n-1}  \xi_\varphi^{n-1}\| \right)
\end{aligned}
\end{equation}
with $z_n$ and $y_n$ defined in \eqref{def_zn_yn}.
\end{lemma}
\proof
According to  \eqref{term_twn}, the term $\mathcal{J}_m^{W}$ may
be bounded as follows
\begin{equation}
\begin{aligned}
\mathcal{J}_m^{W} \leq &  \sum_{n=1}^{m} \frac{2}{k_n +k_{n-1}} \int_{t^{n-1}}^{t^n} |s-t^{n-\frac 1 2}|\,( |\<z_n, \hat\rho(s)\>|+|\<y_n, \hat\rho(s)\>|)\,ds
 \\& + \sum_{n=1}^{m} \int_{t^{n-1}}^{t^n} |s-t^{n-\frac 1 2}|\, ( |\< \xi_\varTheta^{n}, \hat\rho(s)\>|+ |\<\pi^n  \xi_\varTheta^{n-1}, \hat\rho(s)\>|) \, ds
 \\& + \sum_{n=1}^{m} \int_{t^{n-1}}^{t^n} |s-t^{n-\frac 1 2}|\, (|\< P_0^n\xi_\varphi^{n}, \hat\rho(s)\>|+ |\<\pi^n  P_0^{n-1} \xi_\varphi^{n-1}, \hat\rho(s)\>|)\,ds ,
\end{aligned}
\end{equation}
and the claimed result follows.
\qed

We can thus conclude the following a posteriori estimates for the parabolic error:

\begin{lemma}[An $L^\infty(L^2)$ a posteriori error bound for $\hat\rho$ - two time-level reconstruction]
\label{final_est_rho_1} For 
$m = 1,\ldots, N,$ the following estimate holds
\begin{equation}
\begin{aligned}
\max_{t\in[0,t^m]} \|\hat\rho(t)\|+& \left(\int_0^{t^m} |\hat\rho(s)|_1^2\, ds\right)^{1/2} \leq   \sqrt{2}\, \|\hat\rho(0)\|+
\E_m^{T,1}(w_n) 
\\&+ \left \{\bigl(\E_m^{T,2} + \E_m^{S,1}(w_n) + \E_m^{S,2} +\E_m^C+\E_m^{D,1}\bigr)^2 + \bigl(\mathcal{E}_{m}^{D,2}\bigr)^2 \right\}^{1/2}\,,
\end{aligned}
\end{equation}
where
\begin{equation}
 \E_m^{D,2} := \sum_{n=1}^{m}  k_n^{1/2} \zeta_{n,2}
\end{equation}
and $\zeta_{n,2}$ defined in \eqref{zeta_n}.
\end{lemma}
\proof
In view of Theorem \ref{Theorem3_1}, we can easily show that
\begin{equation}
\begin{aligned}
\|\hat \rho(t^m_{\star})\|^2 + \int_0^{t^m} |\hat\rho(s)|_1^2\, ds \leq 2\, \|\hat\rho(0)\|^2 + 2\, \mathcal{J}_m .
\end{aligned}
\end{equation}\bigskip
Thus, by making use of the previous lemmas, we can conclude that
\begin{equation}
\begin{aligned}
\|\hat \rho(t^m_{\star})\|^2 + \int_0^{t^m} &|\hat\rho(s)|_1^2\, ds \leq 2\, \|\hat\rho(0)\|^2 + 2\, \sum_{n=1}^{m} k_n \,\gamma_n^2(w_n) \\
& + 4 \|\hat \rho(t^m_{\star})\|\, \sum_{n=1}^{m} k_n\,  \bigl(\|\xi_{\Theta}^n \|+ k_n \,\eta_n (w_n) + \delta_n  + \beta_n +\zeta_{n,1} + \|\xi_\varphi^n\| \bigr)
\\&+4\sum_{n=1}^{m}  \left(\int_{t^{n-1}}^{t^n}|\hat\rho(s) |_1^2\right)^{1/2}k_n^{1/2} \zeta_{n,2}.
\end{aligned}
\end{equation} 
The final estimate is derived by using the following fact: Let $c\in \Re$ and ${\bf a}=(a_0, a_1, \ldots, a_m), \,{\bf b}=(b_0, b_1, 
\ldots, b_m)\in \Re^{m+1}$ be such that $|{\bf a}|^2 \leq c^2 + {\bf a}
\cdot {\bf b}$, then $|{\bf a}| \leq |c| +  |{\bf b}|.$
Indeed, we apply the above result to the case
\begin{equation*}
\begin{aligned}
c =& \left(2\, \|\hat\rho(0)\|^2 + 2\sum_{n=1}^{m} k_n \gamma_n^2(w_n)\right)^{1/2},\\
a_0=&\,  \|\hat\rho(t^m_\star)\|, \quad a_n= \left (\int_{t^{n-1}}^{t^n} |\hat\rho(s)|_1^2\, ds \right)^{1/2}, \; n=1,\ldots, m,\\
 b_0=&\, 4\sum_{n=1}^{m}  k_n \bigl(\|\xi_{\Theta}^n \|+ k_n \,\eta_n (w_n) + \delta_n  + \beta_n +\zeta_{n,1} + \|\xi_\varphi^n\| \bigr)
\\
b_n=&\,4\, k_n^{1/2} \zeta_{n,2}, \; n=1,\ldots, m,
\end{aligned}
\end{equation*}
to get the final estimate.
\qed

\begin{lemma}[An $L^\infty(L^2)$ a posteriori error bound for $\hat\rho$ - three time-level reconstruction]
\label{final_est_rho_2} For 
$m = 1,\ldots, N,$ the following estimate holds
\begin{equation}
\begin{aligned}
\max_{t\in[0,t^m]} \|\hat\rho(t)\|+& \left(\int_0^{t^m} |\hat\rho(s)|_1^2\, ds\right)^{1/2} \leq   \sqrt{2}\, \|\hat\rho(0)\|+
\E_m^{T,1}(\widetilde{w}_n) + \Bigl \{\bigl(\E_m^{T,2} +\E_m^{T,3}\\&+ \E_m^{S,1}(\widetilde{w}_n) + \E_m^{S,2} +\E_m^C+\E_m^{D,1} + \E_{m,1}\bigr)^2  + \bigl(\mathcal{E}_{m}^{D,2}\bigr)^2 \Bigr\}^{1/2},
\end{aligned}
\end{equation}
where
\begin{equation}
\begin{aligned}
\mathcal{E}_m^{T,3}:= &\sum_{n=1}^{m}   \frac{k_n^2}{2(k_n +k_{n-1})} \,  \|z_n\| \\
 \mathcal{E}_{m,1}:= &\sum_{n=1}^{m}   k_n\,  \Bigl\{ \frac{k_n}{2(k_n +k_{n-1})} \, \|y_n\|+ \frac{k_n}{4}\left(\|\xi_\Theta^n\| + \|\pi^n  \xi_\varTheta^{n-1}\|\right) 
\\ &\qquad+\frac{ k_n}{4} \,  \left( \|P_0^n\xi_\varphi^n\|
 +\ \|\pi^n P_0^{n-1}  \xi_\varphi^{n-1}\| \right)\Bigr\}.  
\end{aligned}
\end{equation}
\end{lemma}

The main result of this paragraph is stated in the next two theorems. 
\bigskip
\begin{theorem}[Final $L^\infty(L^2)$ a posteriori error estimate based on one time subinterval]
\label{Th_L2_err_est} For 
$m = 1,\ldots, N,$ the following estimate holds
\begin{equation}
\begin{aligned}
\max_{t\in[0,t^m]} \|u(t)-U(t)\| \leq &  
 \sqrt{2} \,\|u^0 - \rec^0u^0\| +
\E_m^{T,1}(w_n) 
+\Bigl\{\bigl(\E_m^{T,2} + \E_m^{S,1}(w_n) \\ &+ \E_m^{S,2} +\E_m^C+\E_m^{D,1}\bigr)^2 + \bigl(\mathcal{E}_{m}^{D,2}\bigr)^2 \Bigr\}^{1/2}
  + \E_m^{\emph{rec}}(w_n) 
+ \E_m^{\emph{ell}}\, .
\end{aligned}
\end{equation}
\end{theorem}
\bigskip

\begin{theorem}[$L^\infty(L^2)$ a posteriori error estimate based on two adjacent time intervals]
\label{Th_L2_err_est} For 
$m = 1,\ldots, N,$ the following estimate holds
\begin{equation}
\begin{aligned}
\max_{t\in[0,t^m]} \|u(t)-U(t)\| \leq &  
 \sqrt{2} \,\|u^0 - \rec^0u^0\| +
\E_m^{T,1}(\widetilde{w}_n) 
+\Bigl \{\bigl(\E_m^{T,2} + \E_m^{T,3}+\E_m^{S,1}(\widetilde{w}_n) \\ &+ \E_m^{S,2} +\E_m^C+\E_m^{D,1} + \E_{m,1}\bigr)^2  + \bigl(\mathcal{E}_{m}^{D,2}\bigr)^2 \Bigr\}^{1/2}
+ \E_m^{\emph{rec}}(\widetilde{w}_n) 
+ \E_m^{\emph{ell}}\, .
\end{aligned}
\end{equation}
\end{theorem}

\section{Asymptotic behavior of the estimators}
\label{Sec:numresults}

In this section we study  the asymptotic behavior of the error estimators  and compare this behavior with the true error. For the implementation of the estimators we used the adaptive finite element library ALBERTA-FEM \cite{SS2005}. 

For our purpose,  we consider the heat equation $\eqref{heat_eq}$ on the unit square, $\varOmega = [0,1]^2,$ $T=1$ and the exact solution $u$ be one of the following: 
\begin{itemize}
\item case (1): $u(x,y,t) = \sin(\pi t)\sin(\pi x)\sin(\pi y),$
\item case (2): $u(x,y,t) = \sin(15\pi t)\sin(\pi x)\sin(\pi y) \quad
\text{(fast in time)}, $
\item case (3):
$u(x,y,t) = \sin(0.5\pi t)\sin(10\pi x)\sin(10\pi y) \quad
\text{(fast in space)}.$
\end{itemize}

We take zero initial condition,  $u^0 = 0,$ and  calculate the right-hand side $f$  by applying
the PDE to $u.$ 

We conduct tests on uniform meshes with uniform time steps.  For  the discretization in space we use linear  Lagrange elements.   The computed quantities are: the  {\it error in  the discrete $L^\infty(0,t^m;L^2(\varOmega))$  norm}  
\begin{equation*}
\max_{0\leq n \leq m} \|e^n\| = \max_{0\leq n \leq m} \|u(t^n)-U^n\|,
\end{equation*} 
the {\it total error,} which is dominated by the discrete $L^2(0,t^m;H^1(\varOmega))$ error, 
\begin{equation*}
e_{\rm total}(t^m) := \max_{0\leq n \leq m} \left( \|e^n\|^2  + \sum_{n=1}^m k_n
\|\nabla e^n\|^2\right)^{1/2},
\end{equation*} 
and  almost all the estimators introduced in Section \ref{sec:errest}. 

\begin{figure}[h!tb]
\begin{center}
\subfigure[{\scriptsize Problem case (2)}]{\includegraphics[scale=0.37]{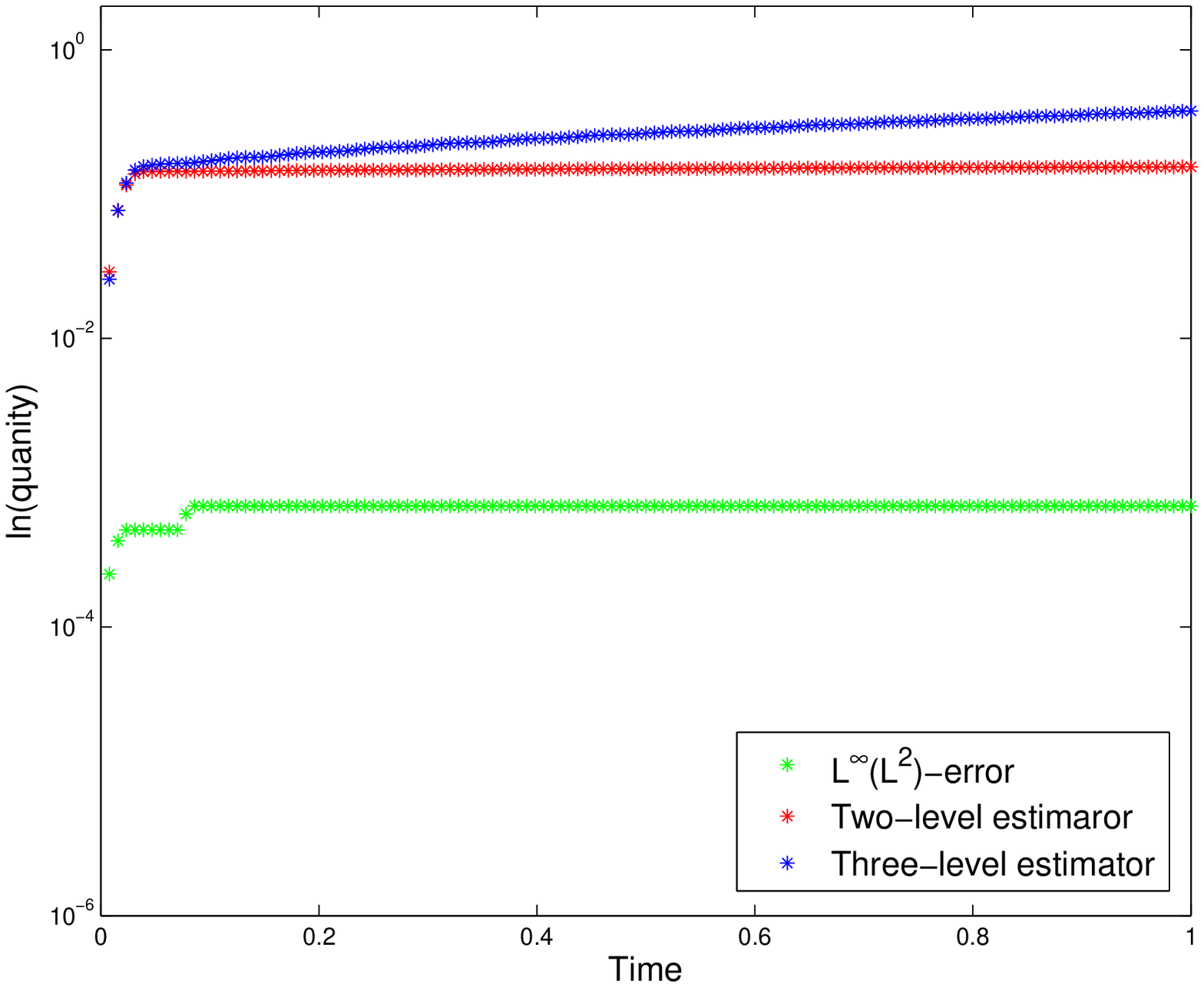}}
\subfigure[{\scriptsize Problem case (3)}]{\includegraphics[scale=0.37]{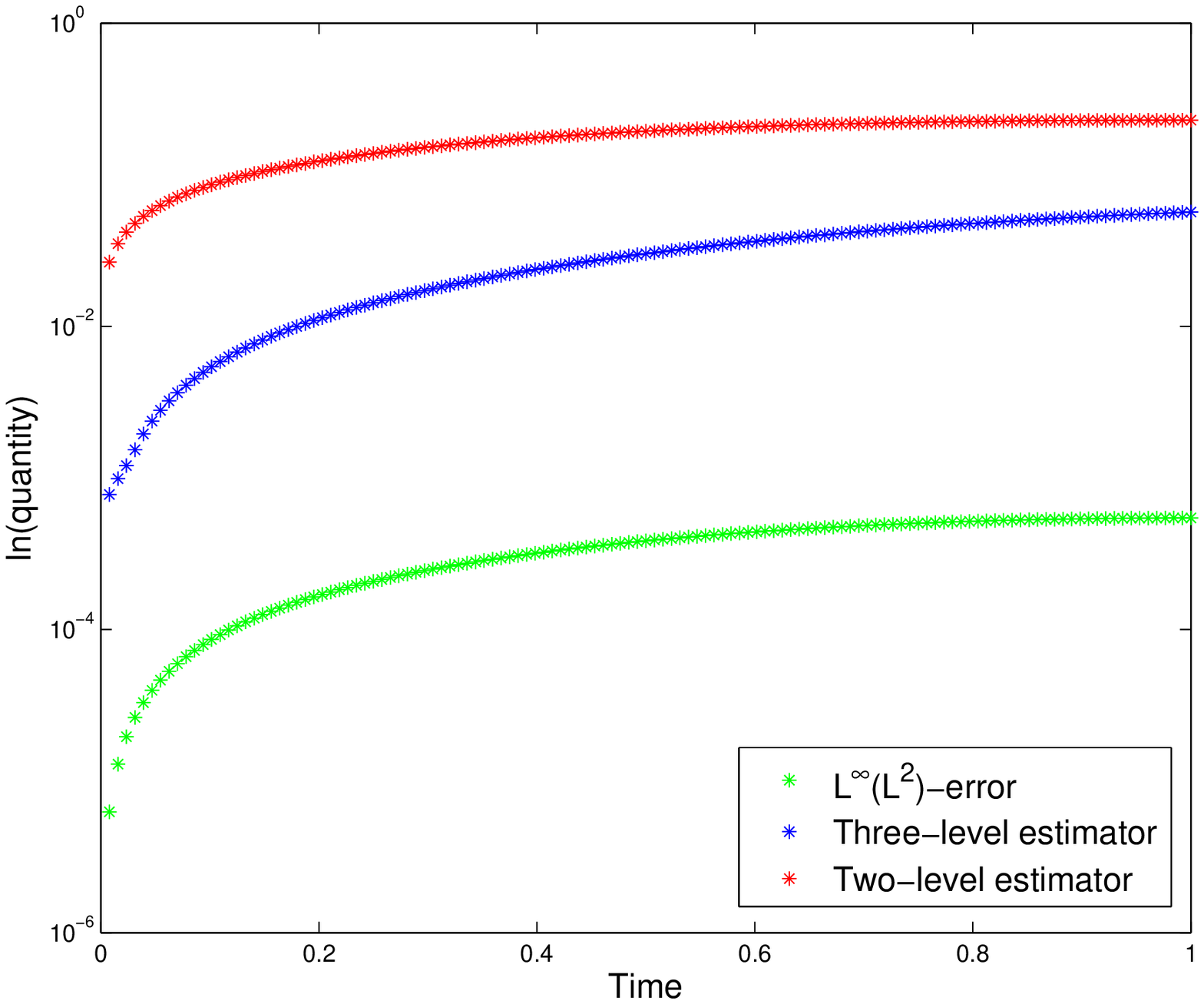}}
\end{center}
\label{est_p2_p3}
\caption{The $L^\infty(0,t^m, L^2(\varOmega))$-error, the two-level estimator $  \E_m$ and the three-level estimator $ \widetilde{\E}_m.$}

\end{figure}

We exclude from the numerical experiments the coarsening error estimator $\E_N^C$ that vanishes as well as  the terms corresponding  to the approximation of data $u^0$ and $f$  which clearly are of optimal order and thus do not contain  interesting information for our purposes.
Moreover, we compute the {\emph{total estimators}} $ \E_{m}$ and $\widetilde\E_{m}, $ which correspond to two time-level reconstruction and three time-level reconstruction, respectively,  defined as follows
\begin{equation}
 \E_{m} := \E_m^{T,1} (w_n)+  \E_m^{T,2} +  \E_m^{S,1}(w_n) +  \E_m^{S,2}  + \E_m^{\text{ell}} + \E_m^{\text{rec}} (w_n), \quad 1\leq m \leq N,
 \end{equation}
and
$$ \widetilde\E_{m} := \E_m^{T,1} (\widetilde{w}_n)+  \E_m^{T,2} +  \E_m^{T,3} +  \E_m^{S,1}(\widetilde{w}_n) +  \E_m^{S,2}  + \E_m^{\text{ell}} + \E_m^{\text{rec}} (\widetilde{w}_n), \quad 1\leq m \leq N.$$
Then, the corresponding effectivity indices are defined as
$$ EI(t^m) := \frac {\E_m}{\max_{0\leq n \leq m} \|e^n\|} \quad \text{and} \quad  \tilde{EI}(t^m) := \frac {\tilde\E_m}{\max_{0\leq n \leq m} \|e^n\|}\,, \quad 1\leq m \leq N .$$

For all quantities of interest    we look at their experimental order of convergence (EOC). The EOC is defined as follows: for a given finite sequence of successive runs (indexed by $i$), the EOC of the  corresponding sequence of  quantities of interest $E(i)$ (error, estimator or part of an estimator), is itself a sequence defined by
$${\rm EOC}(E(i))  = \frac {\log(E(i + 1)/E(i))} { \log(h(i + 1)/h(i)) }\, ,$$
where $h(i)$ denotes the mesh-size of the run $i.$ The values of EOC of a given quantity indicates its order. 
\begin{figure}[h!tb]
\begin{center}
\subfigure[{\scriptsize Problem case (2)  }]{\includegraphics[scale=0.37]{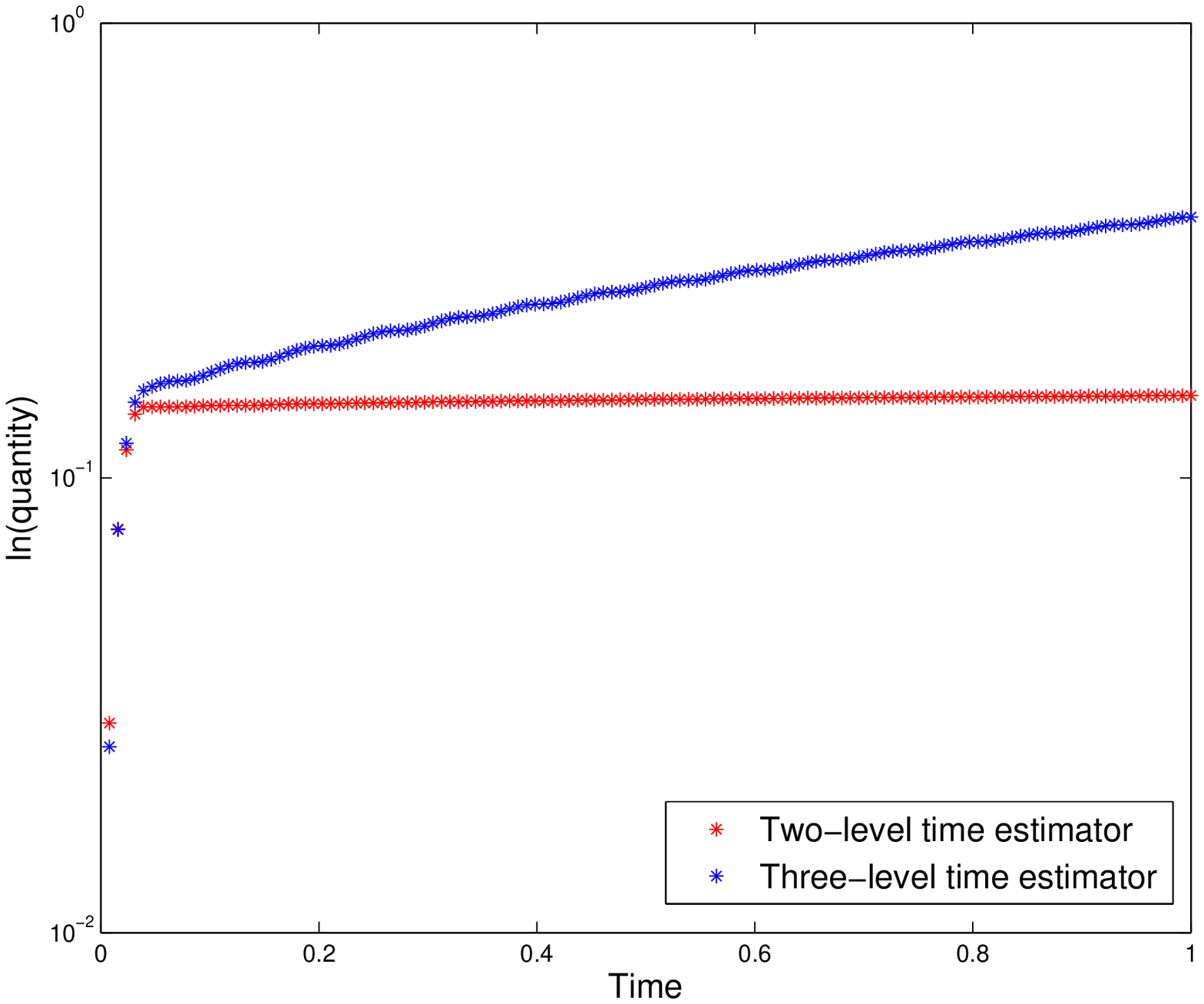}}
\subfigure[{\scriptsize Problem case (3)}]{\includegraphics[scale=0.37]{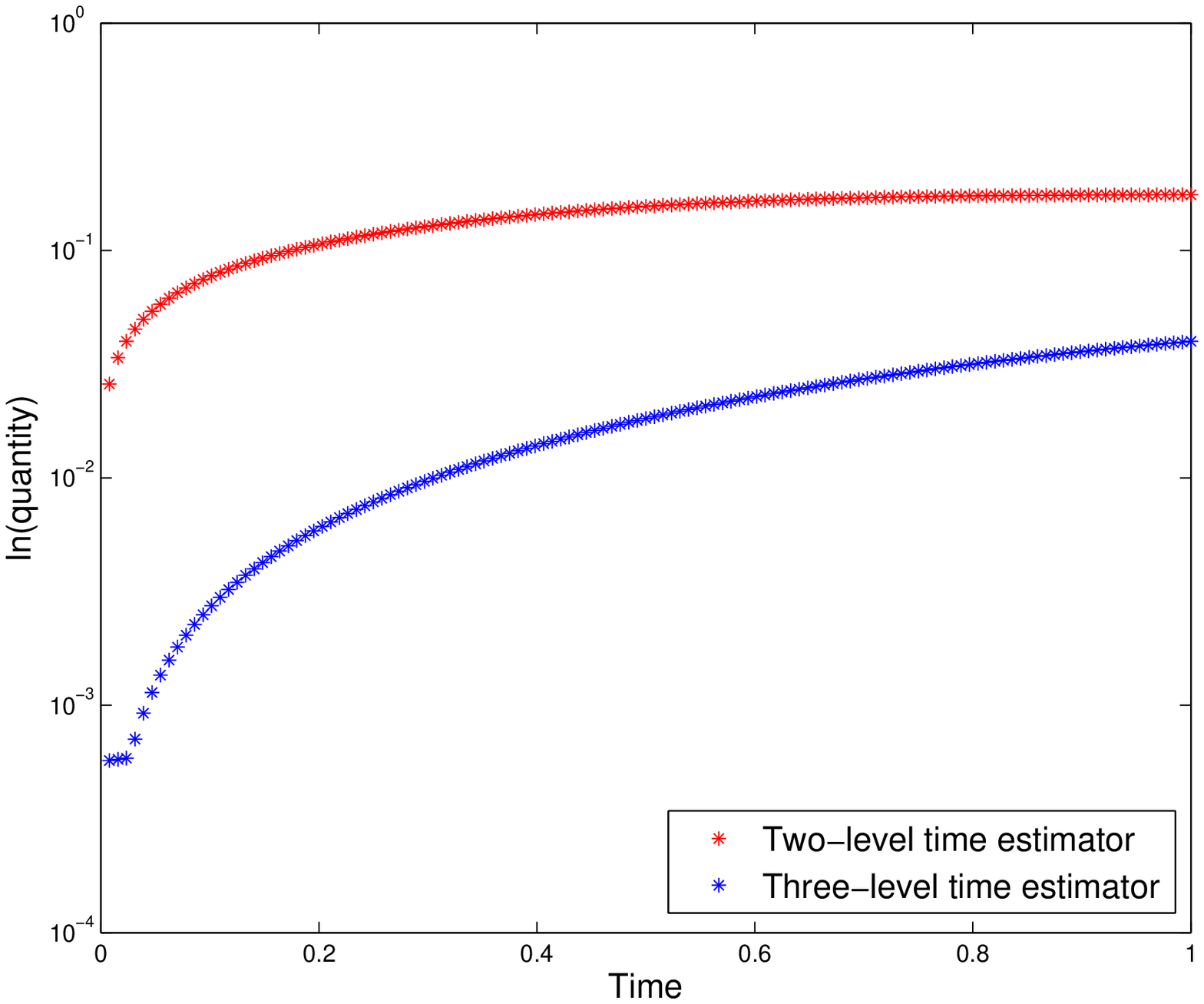}}
\end{center}
\label{time_est_p2_p3}
\caption{The two-level time estimator $ \E_m^{T,1} (w_n)+  \E_m^{T,2}+ \E_m^{\text{rec}} (w_n)$ and the three-level time estimator $ \E_m^{T,1} (\widetilde{w}_n)+  \E_m^{T,2} +  \E_m^{T,3}+ \E_m^{\text{rec}} (\widetilde{w}_n).$}
\end{figure}

\begin{figure}[h!tb]
\begin{center}
\subfigure[{\scriptsize Problem case (2) }]{\includegraphics[scale=0.37]{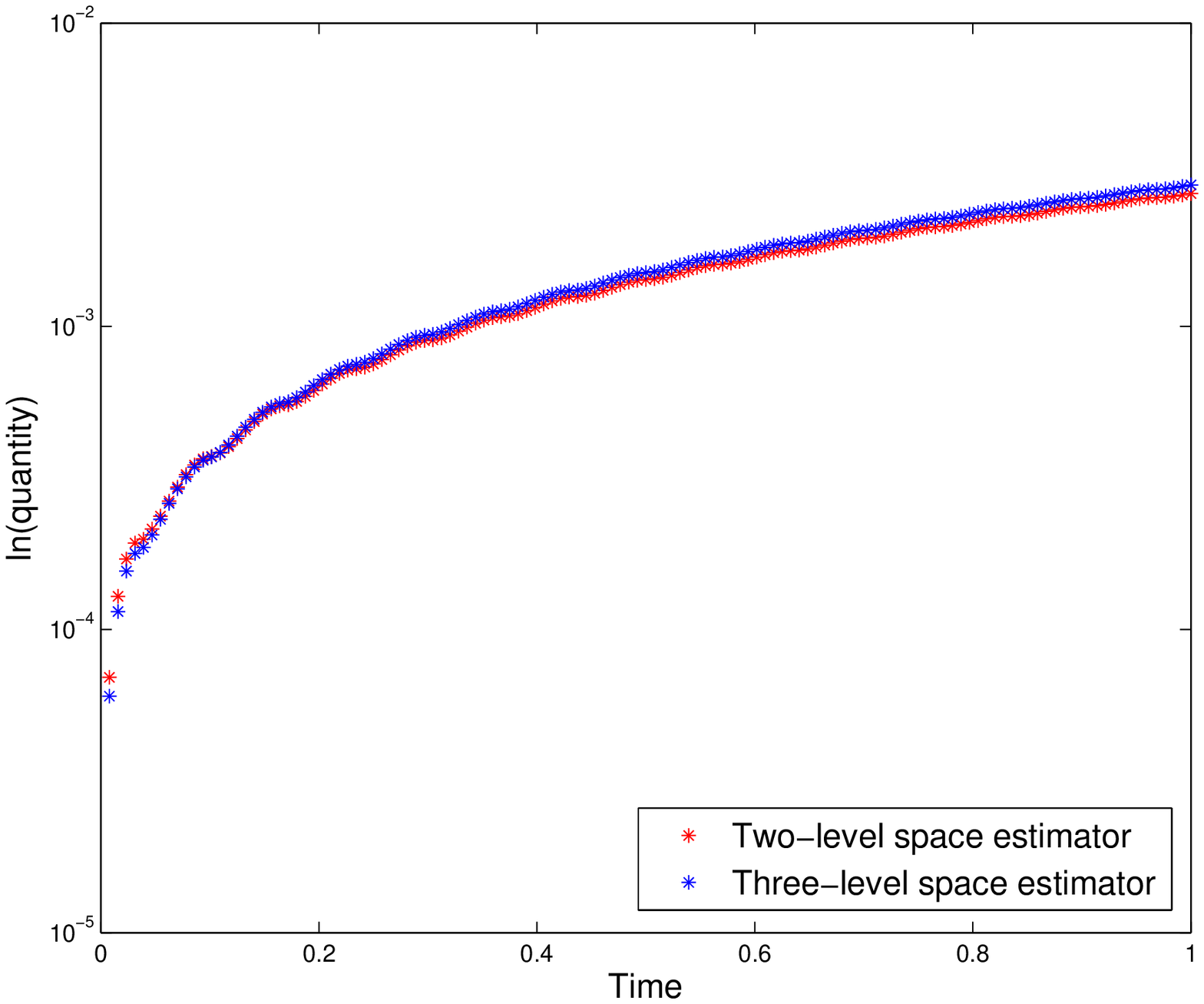}}
\subfigure[{\scriptsize Problem case (3) }]{\includegraphics[scale=0.37]{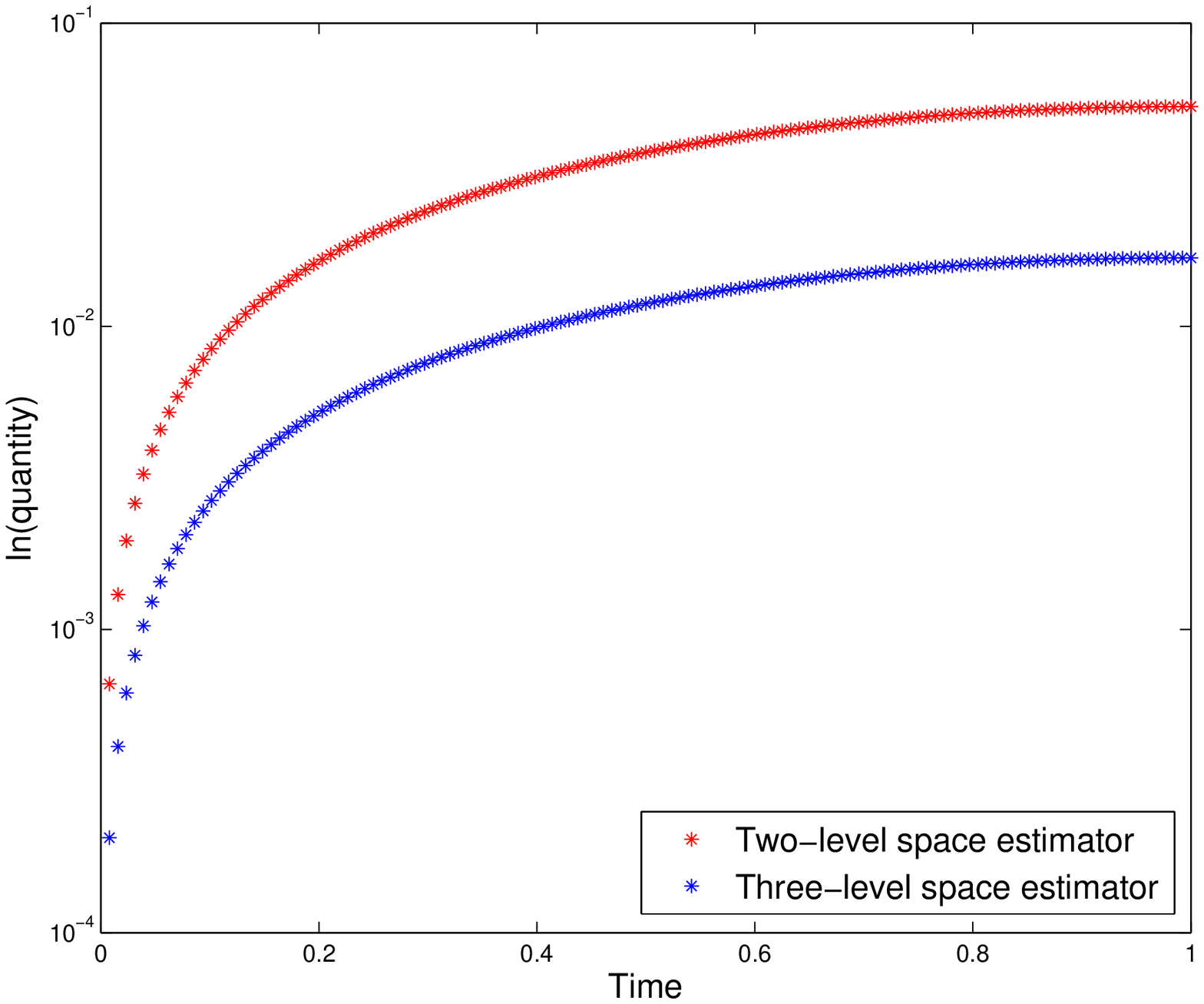}}
\end{center}
\label{space_est_p2_p3}
\caption{The two-level space estimator $  \E_m^{S,1}(w_n) +  \E_m^{S,2}  + \E_m^{\text{ell}} $ and the three-level space estimator $ \E_m^{S,1}(\widetilde{w}_n) +  \E_m^{S,2}  + \E_m^{\text{ell}}.$}
\end{figure}

\begin{table}[!ht]
\begin{tabular}{c}
Errors\\
\end{tabular}
\\
\begin{tabular}{||c||c|c||c|c||c|c||c|c||}
\hline
h = k & $\max_n \|e^n\|$ & EOC & $e_{\rm total}(t^N)$ & EOC &$\E_N$ &$ EI(t^N)$  &$\widetilde{\E}_N$ & $\widetilde{EI}(t^N)$\\
\hline \hline
1.2500e-01  & 1.4481e-03 &         & 3.7925e-02 &         &3.6810e-01 & 254 &4.9986e-01&345\\
6.2500e-02  & 3.4561e-04 & 2.04 & 1.8534e-02 & 1.03 &8.6375e-02 & 249 &1.2473e-01&360\\
3.1250e-02  & 8.4256e-05 & 2.02 & 9.1731e-03 & 1.01& 2.0828e-02 & 247 &3.1030e-02&368\\
1.5625e-02  & 2.0821e-05 & 2.01 & 4.5645e-03 & 1.00& 4.7780e-03 & 229 &7.7309e-03&371\\
7.8125e-03  & 5.1718e-06 & 2.00 & 2.2769e-03 & 1.00& 1.1813e-03 & 228 &1.9289e-03&372\\
\hline
\end{tabular}
\caption{Problem case (1): the $L^\infty(0,1;L^2(\varOmega))$-error, the total error $e_{\rm total}(t^N),$  and the corresponding EOCs, the two- and three-level estimators $\E_N$ and $\widetilde{\E}_N,$ respectively, and the corresponding effectivity indices $EI(t^N)$ and $\widetilde{EI}(t^N)$. }
\label{rec_est_table}
\end{table}

\begin{table}[!ht]
\begin{tabular}{c}
Reconstruction Error Estimators\\
\end{tabular}
\\
\begin{tabular}{|c||c|c||c|c||c|c||c|c|}
\hline
h = k & $\mathcal{E}_N^{\emph{ell}}$ & EOC & $\mathcal{E}_N^{\emph{rec}}(w_n)$ & EOC &$\mathcal{E}_N^{\emph{rec}}(\widetilde{w}_n)$ & EOC \\
\hline \hline
1.2500e-01  & 2.2702e-02 &         & 3.2184e-02 &         & 9.9205e-03 &\\
6.2500e-02  & 5.6786e-03 & 1.99 & 7.7765e-03 & 2.03 & 2.4275e-03& 2.02\\
3.1250e-02  & 1.4198e-03 & 2.00 & 1.9227e-03 & 2.01 & 6.0352e-04& 2.00\\
1.5625e-02  & 3.5497e-04 & 2.00 & 4.7951e-04 & 2.00 & 1.5066e-04 & 2.00\\
7.8125e-03  & 8.8741e-05 & 2.00 & 1.1979e-04 & 2.00 & 3.7654e-05 & 2.00\\
\hline
\end{tabular}
\caption{Problem case (1): the elliptic reconstruction estimator $\mathcal{E}_N^{\text{ell}},$ the time reconstruction estimators $\mathcal{E}_N^{\text{rec}}(w_n),\mathcal{E}_N^{\text{rec}}(\widetilde{w}_n)$ and the corresponding EOCs. }
\label{rec_est_table}
\end{table}
Since the finite element spaces consist of linear Lagrange elements and the
Crank--Nicolson method is second-order accurate, the error in 
$L^\infty(0,T;L^2(\varOmega))$ norm is $O(k^2+h^2)$.
The main conclusion of this paragraph is that all  error estimators, in both cases of time reconstruction,  decrease with at least second order with respect to time and spatial variable, Tables 1-4.  
\begin{table}
\begin{tabular}{c}
Time Estimators\\
\end{tabular}
\\
\begin{tabular}{|c||c|c||c|c||c|c||c|c|}
\hline
h = k& $\mathcal{E}_m^{T,1}(w_n)$ & EOC & $\E_m^{T,1}(\widetilde{w}_n)$ & EOC &$\mathcal{E}_m^{T,2}$ & EOC &$\E_m^{T,3}$& EOC\\
\hline \hline
1.2500e-01 & 8.5528e-02 &         & 2.6086e-02 &         &  1.5603e-01 &         & 2.3606e-01 & \\
6.2500e-02 & 1.9581e-02 & 2.09 & 6.0501e-03 & 2.07 &  3.9203e-02 & 1.99 & 6.0174e-02 & 1.98\\
3.1250e-02 & 4.6603e-03 & 2.05 & 1.4486e-03 & 2.04 &  9.7796e-03 & 2.00 & 1.5116e-02 & 1.99\\
1.5625e-02 & 1.1340e-03 & 2.02 & 3.5392e-04 & 2.02 &  2.4396e-03 & 2.00 & 3.7834e-03 & 1.99\\
7.8125e-03 & 2.7935e-04 & 2.01 & 8.7434e-05 & 2.01 &  6.0908e-04 & 2.00 & 9.4614e-04 & 2.00\\
\hline
\end{tabular}
\caption{Problem case (1): the time estimators $\mathcal{E}_N^{T,1}(w_n),  \E_N^{T,1}(\widetilde{w}_n), \mathcal{E}_N^{T,2}, \E_N^{T,3}$ and the corresponding EOCs.   }
\label{time_est_table}
\end{table}

\begin{table}
\begin{tabular}{c}
Space Estimators\\
\end{tabular}
\\
\begin{tabular}{|c||c|c||c|c||c|c|}
\hline
h = k & $\mathcal{E}_N^{S,1}(w_n)$ & EOC & $\E_N^{S,1}(\widetilde{w}_n)$ & EOC &$\mathcal{E}_N^{S,2}$  & EOC \\
\hline \hline
1.2500e-01& 3.1258e-02 &         & 8.6889e-03 &         & 4.0374e-02 &          \\
6.2500e-02& 4.0411e-03 & 2.78 & 1.1079e-03 & 2.80 & 1.0095e-02 & 2.00  \\
3.1250e-02& 5.2213e-04 & 2.78 & 1.3917e-04 & 2.82 & 2.5239e-03 & 2.00  \\
1.5625e-02& 6.7829e-05 & 2.77 & 1.7417e-05 & 2.82 & 6.3099e-04 & 2.00  \\
7.8125e-03& 8.7779e-06 & 2.77 & 2.1778e-06 & 2.82 & 1.5775e-04 & 2.00  \\
\hline
\end{tabular}
\caption{Problem case (1): the space estimators $\mathcal{E}_N^{S,1}(w_n),  \mathcal{E}_N^{S,1}(\widetilde{w}_n), \mathcal{E}_N^{S,2}$ and the corresponding EOCs. }
\label{space_est_table}
\end{table}

\bibliographystyle{abbrv}
\bibliography{bibfile_proposal,bibfile_stokes,bibfile_theta}

\begin{thebibliography}{10}

\bibitem{AMN2006}
G.~Akrivis, C.~Makridakis, and R.~H. Nochetto.
\newblock A posteriori error estimates for the {C}rank--{N}icolson method for
  parabolic equations.
\newblock {\em Math. Comp.}, 75:511--531, 2006.

\bibitem{BKM2012}
E.~B\"{a}nsch, F.~Karakatsani, and C.~Makridakis.
\newblock A posteriori error control for fully discrete {C}rank--{N}icolson
  schemes.
\newblock {\em SIAM J. Numer. Anal.}, 6:2845--2872, 2012.

\bibitem{BKM2013}
E.~B\"{a}nsch, F.~Karakatsani, and C.~Makridakis.
\newblock The effect of mesh modification in time on the error control of fully
  discrete approximations for parabolic equations.
\newblock {\em Appl. Numer. Math.}, 67:35--63, 2013.

\bibitem{Bri-Glo_Per87}
M.~Bristeau, R.~Glowinski, and J.~Periaux.
\newblock Numerical methods for {N}avier--{S}tokes equations. applications to
  the simulation of compressible and incompressible viscous flows.
\newblock {\em Computer Physics Reports}, 6:73--187, 1987.

\bibitem{EJ_l1}
K.~Eriksson and C.~Johnson.
\newblock Adaptive finite element methods for parabolic problems. {I}. {A}
  linear model problem.
\newblock {\em SIAM J. Numer. Anal.}, 28(1):43--77, 1991.

\bibitem{EJ_nonl}
K.~Eriksson and C.~Johnson.
\newblock Adaptive finite element methods for parabolic problems. {I}{V}.
  {N}onlinear problems.
\newblock {\em SIAM J. Numer. Anal.}, 32(6):1729--1749, 1995.

\bibitem{Glowinski85}
R.~Glowinski.
\newblock Viscous flow simulations by finite element methods and related
  numerical techniques.
\newblock In E.~Murman and S.~Abarbanel, editors, {\em Progress in
  Supercomputing in Computational Fluid Dynamics}, pages 173--210.
  Birkh\"auser, Boston, 1985.

\bibitem{Glowinski86}
R.~Glowinski.
\newblock Splitting methods for the numerical solution of the incompressible
  {N}avier--{S}tokes equations.
\newblock In A.~D. A.V.~Balakrishman and J.~Lions, editors, {\em Vistas in
  Applied Mathematics}, pages 57--95. Optimization Software, New York, 1986.

\bibitem{Glowinski:2003}
R.~Glowinski.
\newblock Finite element methods for incompressible viscous flow.
\newblock In {\em Handbook of numerical analysis, Vol. IX}. North-Holland,
  Amsterdam, 2003.

\bibitem{JNT}
C.~Johnson, Y.~Y. Nie, and V.~Thom{\'e}e.
\newblock An a posteriori error estimate and adaptive timestep control for a
  backward {E}uler discretization of a parabolic problem.
\newblock {\em SIAM J. Numer. Anal.}, 27(2):277--291, 1990.

\bibitem{K2012}
F.~Karakatsani.
\newblock A posteriori error estimates for the fractional-step
  $\vartheta$-scheme for linear parabolic equations.
\newblock {\em IMA J. Numer. Anal.}, 32(1):141--162, 2012.

\bibitem{LM2006}
O.~Lakkis and C.~Makridakis.
\newblock Elliptic reconstruction and a posteriori error estimates for fully
  discrete linear parabolic problems.
\newblock {\em Math. Comp.}, 75(256):1627--1658, 2006.

\bibitem{LPP2009}
A.~Lozinski, M.~Picasso, and V.~Prachittham.
\newblock An anisotropic error estimator for the {C}rank-{N}icolson scheme.
\newblock {\em SIAM J. Sci. Comp}, 31(4):2757--2783, 2009.

\bibitem{MN2003}
C.~Makridakis and R.~H. Nochetto.
\newblock Elliptic reconstruction and a posteriori error estimates for
  parabolic problems.
\newblock {\em SIAM J. Numer. Anal.}, 41(4):1585--1594, 2003.

\bibitem{MR2014}
D.~Meidner and T.~Richter.
\newblock Goal-oriented error estimation for the fractional step theta scheme.
\newblock {\em Computational Methods in Applied Mathematics}, 2014.

\bibitem{Mueller94}
S.~M\"uller, A.~Prohl, R.~Rannacher, and S.~Turek.
\newblock Implicit time-discretization of the nonstationary incompressible
  navier-stokes equations.
\newblock In W.~Hackbusch and G.~Wittum, editors, {\em Proc. Workshop. ``Fast
  Solvers for Flow Problems''}, pages 175--191, Kiel, Germany, Jan. 14-16 1994.
  Vieweg, Braunschweig.

\bibitem{Urbaniak93}
S.~M\"uller-Urbaniak.
\newblock {\em Eine Analyse des Zwischenschritt-$\theta$-Verfahrens zur
  L\"osung der instation\"aren Navier-Stokes-Gleichungen}.
\newblock PhD thesis, University of Heidelberg, 1993.

\bibitem{NSaV}
R.~H. Nochetto, G.~Savar{\'e}, and C.~Verdi.
\newblock A posteriori error estimates for variable time-step discretizations
  of nonlinear evolution equations.
\newblock {\em Comm. Pure Appl. Math.}, 53(5):525--589, 2000.

\bibitem{Prachittham2009}
V.~Prachittham.
\newblock {\em Space-time adaptive algorithms for parabolic problems: a
  posteriori error estimates and application to microfluidics}.
\newblock Ph. d. thesis, EPFL, Laussane, 2009.

\bibitem{Rannacher98}
R.~Rannacher.
\newblock Numerical analysis of nonstationary fluid flow (a survey).
\newblock In V.~Boffi and H.~Neunzert, editors, {\em Applications of
  Mathematics in Industry and Technology}, pages 34--53. B.G. Teubner,
  Stuttgart, 1998.

\bibitem{SS2005}
A.~Schmidt and K.~Siebert.
\newblock {\em Design of adaptive finite element software: The finite element
  toolbox {ALBERTA}}, volume~42 of {\em Springer LNCSE Series}.
\newblock Springer-Verlag, Berlin, 2005.

\bibitem{scott-zhang}
L.~R. Scott and S.~Zhang.
\newblock Finite element interpolation of nonsmooth functions satisfying
  boundary conditions.
\newblock {\em Math. Comp.}, 54(190):483--493, 1990.

\bibitem{Thomee97}
V.~Thom{\'e}e.
\newblock {\em Galerkin finite element methods for parabolic problems}.
\newblock Springer-Verlag, Berlin, 1997.

\bibitem{Verfuerth2003}
R.~Verf\"urth.
\newblock A posteriori error estimates for finite element discretizations of
  the heat equation.
\newblock {\em Calcolo}, 40:195--212, 2003.

\bibitem{Wheeler73}
M.~F. Wheeler.
\newblock A priori ${L}\sb{2}$ error estimates for {G}alerkin approximations to
  parabolic partial differential equations.
\newblock {\em SIAM J. Numer. Anal.}, 10:723--759, 1973.

\end{thebibliography}
\end{document}